\documentclass[11pt]{article}

\title{Lecture notes on locally compact quantum groups}
\author{Alfons Van Daele}

\usepackage{amsmath} 
\usepackage{amsthm} 
\usepackage{amssymb} 

\usepackage{hyperref}

\newtheorem{thm}{Theorem}
\newtheorem{inspr}[thm]{}

\newenvironment{definitie}{\begin{itemize}\item[ ]\hspace{-26pt}\bf Definition \rm }{\end{itemize}}
\newenvironment{notatie}{\begin{itemize}\item[ ]\hspace{-26pt}\bf Notation \rm }{\end{itemize}}
\newenvironment{voorbeeld}{\begin{itemize}\item[ ]\hspace{-26pt}\bf Example \rm }{\end{itemize}}
\newenvironment{stelling}{\begin{itemize}\item[ ]\hspace{-26pt}\bf Theorem \rm }{\end{itemize}}
\newenvironment{propositie}{\begin{itemize}\item[ ]\hspace{-26pt}\bf Proposition \rm }{\end{itemize}}
\newenvironment{lemma}{\begin{itemize}\item[ ]\hspace{-26pt}\bf Lemma \rm }{\end{itemize}}
\newenvironment{opmerking}{\begin{itemize}\item[ ]\hspace{-26pt}\bf Remark \rm }{\end{itemize}}
\newenvironment{voorwaarde}{\begin{itemize}\item[ ]\hspace{-26pt}\bf Condition \rm }{\end{itemize}}

\renewcommand{\Bbb}{\mathbb} 

\newcommand{\defin}{\begin{inspr}\begin{definitie}}  
\newcommand{\edefin}{\end{definitie}\end{inspr}}

\newcommand{\notat}{\begin{inspr}\begin{notatie}}  
\newcommand{\enotat}{\end{notatie}\end{inspr}}

\newcommand{\voorb}{\begin{inspr}\begin{voorbeeld}}  
\newcommand{\evoorb}{\end{voorbeeld}\end{inspr}}

\newcommand{\stel}{\begin{inspr}\begin{stelling}}
\newcommand{\estel}{\end{stelling}\end{inspr}}

\newcommand{\prop}{\begin{inspr}\begin{propositie}}
\newcommand{\eprop}{\end{propositie}\end{inspr}}

\newcommand{\lem}{\begin{inspr}\begin{lemma}}
\newcommand{\elem}{\end{lemma}\end{inspr}}

\newcommand{\opm}{\begin{inspr}\begin{opmerking}}
\newcommand{\eopm}{\end{opmerking}\end{inspr}}

\newcommand{\voorw}{\begin{inspr}\begin{voorwaarde}}
\newcommand{\evoorw}{\end{voorwaarde}\end{inspr}}

\newcommand{\bew}{\vspace{-0.3cm}\begin{itemize}\item[ ] \bf Proof\rm: }
\newcommand{\ebew}{\hfill $\qed$ \end{itemize}}

\newcommand{\snl}{\vskip 7pt} 
\newcommand{\nl}{\vskip 12pt}

\newcommand{\ot}{\otimes}
\newcommand{\di}{\diamond}

\newcommand{\tr}{\triangleright}

\newcommand{\tussenen}{\qquad\quad\text{and}\qquad\quad}

\numberwithin{thm}{section}   
\numberwithin{equation}{section} 

\begin{document}

\centerline{\bf \Large Algebraic quantum groupoids - An example}
\vspace{13pt}
\centerline{\it  Alfons Van Daele \rm $^{(*)}$}

\bigskip\bigskip

{\bf Abstract} 
\nl
Let $B$ and $C$ be non-degenerate idempotent algebras and assume that $E$ is a regular separability idempotent in $M(B\ot C)$. Define $A=C\ot B$ and $\Delta:A\to M(A\ot A)$ by $\Delta(c\ot b)=c\ot E\ot b$. The pair $(A,\Delta)$ is a {\it weak multiplier Hopf algebra}. Because we assume that $E$ is regular, it is a {\it regular weak multiplier Hopf algebra}. 
\snl
There is a faithful left integral  on $(A,\Delta)$ that is also right invariant. Therefore, we call $(A,\Delta)$ a {\it unimodular} algebraic quantum groupoid. By the general theory, the dual $(\widehat A,\widehat \Delta)$ can be constructed and it is again an algebraic quantum groupoid.
\snl
In this paper, we treat this algebraic quantum groupoid and its dual in great detail. The main purpose is to illustrate various aspects of the general theory. For this reason, we will also recall the basic notions and results of separability idempotents and weak multiplier Hopf algebras with integrals.
\snl
The paper is to be considered as an expository note.
\nl
Date: {\it 16 February 2017}
\nl

\vskip 8cm
\hrule
\medskip
\begin{itemize}
\item[$^{(*)}$] Department of Mathematics, University of Leuven, Celestijnenlaan 200B,\newline
B-3001 Heverlee, Belgium. E-mail: \texttt{alfons.vandaele@kuleuven.be}
\end{itemize} 

\newpage

\setcounter{section}{-1}  

\section{\hspace{-17pt}. Introduction} \label{s:introduction}  

Let $G$ be a groupoid. We use $s$ for the source map and $t$ for the target map. And we consider the units as sitting inside $G$.
\snl
There are canonically associated two weak multiplier Hopf algebras. For the first one, we take for $A$ the algebra $K(G)$ of complex functions with finite support in $G$ and pointwise operations. The coproduct $\Delta$ on $A$ is defined by
\begin{equation*}
\Delta(f)(p,q)=
\begin{cases}
		f(pq) & \text{if $pq$ is defined},\\
		0 & \text{otherwise}.
\end{cases}
\end{equation*}
Remark that the algebra $A$ has no identity, except when $G$ is finite. So only for a finite groupoid, we get a weak Hopf algebra here. For an infinite groupoid, we have a weak {\it multiplier} Hopf algebra.
\snl
For the second one, we take  the groupoid algebra $\mathbb C G$. We denote by $p\mapsto \lambda_p$ the canonical embedding of $G$ in this algebra. Then $\mathbb C G$ is the linear span of these elements $\lambda_p$ with $p\in G$. The coproduct  is now defined by $\Delta(\lambda_p)=\lambda_p\ot\lambda_p$. Remark that the groupoid algebra has an identity if and only if the set of units is finite. So here, we only get a weak Hopf algebra if there are only finitely many units in the groupoid.
\snl
Details about these two examples of weak multiplier Hopf algebras can be found e.g.\ in Example 3.1 in (version 2 of) \cite{VD-W5}). 
\snl
In the two cases integrals exist. For the function algebra, we have that $f\mapsto \sum_p g(p)f(p)$ is a left integral if and only if $g$ is a function on $G$ with the property that $g(p)=g(q)$ if $s(p)=s(q)$, i.e.\ when $p$ and $q$ have the same source. On the other hand, the map $f\mapsto \sum_p h(p)f(p)$ will give a right integral if and only if $h$ is a function on $G$ with the property that $h(p)=h(q)$ when $t(p)=t(q)$, i.e.\ if $p$ and $q$ have the same target.  If we take $f\mapsto \sum_pf(p)$ we get a faithful left integral that is also a right integral. This is proven in Theorem 3.1 of \cite{VD-W6}. 
\snl
For the groupoid algebra, we have that $\lambda_p\mapsto g(p)$ will be a left integral if and only $g$ is a complex function on $G$ with support in the set of units. In this case, the left integrals are the same as the right integrals (as the coproduct is coabelian). Again see Proposition 3.3 of \cite{VD-W6}. 
\snl
In Section 3 of \cite{VD-W6} it is also shown that these two weak multiplier Hopf algebras are dual to each other in the sense of duality for weak multiplier Hopf algebras as studied in \cite{VD-W6}. We used this example to illustrate the construction of the dual that we developed in that paper.
\snl
In this paper, we want to present another example to illustrate this procedure. The example is already described in Section 3 of \cite{VD-W6}, but without the details. This note is devoted to a detailed treatment of this example.
\snl
The example is a generalization (in another direction) of the following special case of the groupoid example above.
\nl
\bf A special case \rm
\nl
Let $X$ be any set. Consider the groupoid $G$ of pairs $(x,y)$ of elements in $X$ with multiplication $pq=(x,y')$, defined for $p=(x,y)$ and $q=(x',y')$ if $y=x'$. The set of units is $X$. The source map $s$ is $(x,y)\mapsto y$ while the target map $t$ is $(x,y)\mapsto x$. The units are considered as sitting in $G$ via the map $x\mapsto (x,x)$. 
\snl
The function algebra $A$ for this groupoid $G$ is the algebra $K(X\times X)$ of complex functions with finite support on $X\times X$ and pointwise product. The algebra $A\ot A$ is identified with $K(X^4)$ where we use $X^4$ for the Cartesian product of 4 copies of $X$. The multiplier algebra $M(A\ot A)$ is the algebra $C(X^4)$ of all complex functions on $X^4$. The coproduct is given by

\begin{equation*}
\Delta(f)(x,y;x',y')=
\begin{cases}
		f(x,y') & \text{if $y=x'$},\\
		0 & \text{otherwise}.
\end{cases}
\end{equation*}

For a left integral $\varphi$ we must have that $(\iota\ot\varphi)\Delta(f)$ is a function in $A_t$ (the multiplier algebra of the target algebra). In this case $A_t$ consist of functions $(x,y)\mapsto h(x)$ for $h\in C(X)$. For a right integral $\psi$, we must have that $(\psi\ot\iota)\Delta(f)$ is a function in $A_s$ (the multiplier algebra of the source algebra). Here $A_s$ is the algebra of functions $(x,y)\mapsto g(y)$ for $g\in C(X)$. Then one can verify that any left integral is of the form $f\mapsto \sum_{x,y} g(y)f(x,y)$ for some $g\in C(X)$  and that any right integral has the form $f\mapsto \sum_{x,y}h(x)f(x,y)$ for some $h\in C(X)$. With $h=1$ and $g=1$ we find the faithful left integral $f\mapsto \sum_{x,y} f(x,y)$, which is also a right integral.
\snl
Next, look at the groupoid algebra $\mathbb C G$ of this (rather trivial) groupoid. In this case, it is the algebra $K(X\times X)$ with a convolution type product
$(fg)(x,y)=\sum f(x,u)g(u,y)$. Now, the coproduct is defined as
\begin{equation*}
\Delta(f)(x,y;x',y')=
\begin{cases}
		f(x,y) & \text{if $x=x'$ and $y=y'$},\\
		0 & \text{otherwise}.
\end{cases}
\end{equation*}

The source and target algebras are given by functions with support in the diagonal. And every integral $\varphi$ is of the form
$f\mapsto \sum_u g(u)f(u,u)$ for some function $g\in C(X)$.
\snl
One can verify that the two weak multiplier Hopf algebras are dual to each other and that the duality is realized with the obvious pairing $(f,g)\mapsto \sum_{x,y} f(x,y)g(x,y)$. This is a consequence of the general result on the duality between $K(G)$ and $\mathbb C G$ for any groupoid, as studied in Section 3 of \cite{VD-W6}. 
\nl
\bf The weak multiplier Hopf algebra of a separability idempotent \rm
\nl
In this paper, we now will study a generalization of the above example. 
\snl
Take two non-degenerate algebras $B$ and $C$ and assume that they are idempotent. We let $E\in M(B\ot C)$ be a regular separability idempotent. We consider the algebra $A=C\ot B$ and define a coproduct $\Delta$ on $A$ by 
\begin{equation*}
\Delta(c\ot b)=c\ot E\ot b
\end{equation*}
for $b\in B$ and $c\in C$. Then $(A,\Delta)$ is a regular weak multiplier Hopf algebra. We will recall this in the preliminary section, Section \ref{s:preliminaries}, of this paper. 
\snl
Let $B=K(X)$ and $C=K(X)$ for a given set $X$. Let $E$ be defined by $E(x,y)=0$ for all $(x,y)$, except when $x=y$ when it is defined $E(x,x)=1$. This is a regular separability idempotent and the weak multiplier Hopf algebra $(A,\Delta)$ is precisely the one defined  in the previous item.
\snl
The weak multiplier Hopf algebra $(A,\Delta)$ associated with a regular separability idempotent $E$ is mentioned already in Section 4 of \cite{VD-W5}. In this paper however, we give a detailed coverage of this example. 
\snl
We discuss the counit, the antipode and the source and target maps and source and target algebras in {\it Section} \ref{s:thewmha} in detail and we give a characterization of left and right integrals. In particular, we show that there is a faithful left integral that is also right invariant.
\snl
In {\it Section} \ref{s:duality} we construct the dual by applying the general theory developed in Section 2 of \cite{VD-W6}. Also here, we give a detailed treatment as it is our intention to illustrate the various steps of the general theory.
\snl
We conclude this note with {\it Section} \ref{s:conclusions} where  we also discuss some further possible research. 
\snl
Note that the whole paper is intended to illustrate the duality procedure as given in \cite{VD-W6} for regular weak multiplier Hopf algebras in the special case of an example. By treating it in detail, the paper becomes partly of an expository nature. 
\nl
\bf Conventions and notations \rm
\snl
We only work with algebras $A$ over $\Bbb C$ (although we believe that this is not essential). We do not assume that they are unital but we need that the product is non-degenerate. We also assume our algebras to be idempotent. This means that every element in the algebra is a sum of products of two elements. The condition is written as $A^2=A$. In fact, the algebras we deal with have local units. Then of course, the product is automatically non-degenerate and also the algebra is idempotent.
\snl
When $A$ is such an algebra, we use $M(A)$ for the multiplier algebra of $A$. When $m$ is in $M(A)$, then by definition we can define $am$ and $mb$ in $A$ for all $a,b\in A$ and we have $(am)b=a(mb)$. The algebra $A$ sits in $M(A)$ as an essential two-sided ideal and $M(A)$ is the largest algebra with identity having this property. 
\snl
We consider $A\ot A$, the tensor product of $A$ with itself. It is again an idempotent, non-degenerate algebra and we can consider the multiplier algebra $M(A\ot A)$. The same is true for a multiple tensor product. 
We use $\zeta$ for the flip map on $A\ot A$, as well as for its natural extension to $M(A\ot A)$.
\snl
We use $1$ for the identity in any of these  multiplier algebras. On the other hand, we mostly use $\iota$ for the identity map on $A$ (or other spaces), although sometimes, we also write $1$ for  this map. The identity element in a group is denoted by $e$. If $G$ is a groupoid, we will also use $e$ for units. Units are considered as being elements of the groupoid and we use $s$ and $t$ for the source and target maps from $G$ to the set of units. 
\snl
A linear functional $\omega$ on an algebra $A$ is called {\it faithful} if given $a\in A$ we have $a=0$ if either $\omega(ab)=0$ for all $b\in A$ or $\omega(ba)=0$ for all $b\in A$. 
\snl
When $A$ is an algebra, we denote by $A^{\text{op}}$ the algebra obtained from $A$ by reversing the product. When $\Delta$ is a coproduct on $A$, we denote by $\Delta^{\text{cop}}$ the coproduct on $A$ obtained by composing $\Delta$ with the flip map $\zeta$.
\snl
For a coproduct $\Delta$, as we define it in Definition 1.1 of \cite{VD-W4}, we assume that $\Delta(a)(1\ot b)$ and $(a\ot 1)\Delta(b)$ are in $A\ot A$ for all $a,b\in A$. This allows us to make use of the {\it Sweedler notation} for the coproduct. The reader who wants to have a deeper understanding of this, is referred to \cite{VD7} where the use of the Sweedler notation for coproducts that do not map into the tensor product, but rather in its multiplier algebra is explained in detail.

\newpage
\bf Basic references \rm
\nl
For the theory of Hopf algebras, we refer to the standard works of Abe \cite{A} and Sweedler\cite{S}. For multiplier Hopf algebras and integrals on multiplier Hopf algebras, we refer to \cite{VD1} and \cite{VD2}. Weak  Hopf algebras have been studied in \cite{B-N-S} and \cite{B-S} and more results are found in \cite{N} and \cite{N-V1}. Various other references on the subject can be found in \cite{Va}. In particular, we refer to \cite{N-V2} because we will use notations and conventions from this paper when dealing with weak Hopf algebras.
\snl
Weak multiplier Hopf algebras have first been introduced in a preliminary paper \cite{VD-W3}. Then the study has been initiated in \cite{VD-W4}. The source and target maps on the one hand, and the source and target algebras on the other hand, for these weak multiplier Hopf algebras are treated in detail in version 1 of \cite{VD-W5} and again, but in greater generality, in version 2 of \cite{VD-W5}. Finally, integrals and duality of weak multiplier Hopf algebras with integrals are found in \cite{VD-W6}. As we mentioned before, this note is meant as a detailed illustration of some of the phenomena, studied in general in that paper.
\snl
For the theory of groupoids, we refer to \cite{Br}, \cite{H},  \cite{P} and  \cite{R}. 
\nl
\bf Note added while preparing this manuscript \rm
\nl
While writing down this manuscript, we came across the work of Thomas Timmermann, {\it On dualitiy for algebraic quantum groups} \cite{T2} that contains similar results. The relation between his work and ours has to be investigated further. 
\nl\nl
\bf Acknowledgments \rm
\nl
I would like to express my gratitude to S. Wang for motivating me to start the research on weak multiplier Hopf algebras, for the cooperation on the subject and the hospitality during several visits to Nanjing. I am indebted to S.L.\ Woronowicz and P. Hajac for the hospitality during my frequent visits to Warsaw, as well as to my colleagues in Oslo  where part of this work was completed.

\section{\hspace{-17pt}. Preliminaries} \label{s:preliminaries}  

In this section we will recall the basic notions and definitions we need for the rest of the paper. The first concept is that of a separability idempotent. The second one is the notion of a weak multiplier Hopf algebra. We do not give all the results we need, just the most important ones. Results that we will need and are not yet included here, will be recalled as soon as they are used for the first time in the following sections.
\nl
\bf Separability idempotents \rm
\nl
We first recall the main definitions and results from \cite{VD4}. There are two (different) versions of this paper. Except when we explicitly refer to the first version, all internal references are with respect to the second one.
\snl
Assume that $B$ and $C$ are non-degenerate and idempotent algebras. Definition 1.4 from \cite{VD4} says the following. 

\defin 
Let $E$ be an idempotent in the multiplier algebra $M(B\ot C)$ and assume that $E(1\ot b)$ and $(c\ot 1)E$ belong to $B\ot C$ for all $b\in B$ and $c\in C$. Assume also that $E$ is full in the sense that the left leg and the right leg of $E$ are respectively  all of $B$ and all of $C$. Furthermore it is required that there are non-degenerate anti-homomorphisms $S_B:B\to M(C)$ and $S_C:C\to M(B)$ satisfying (and characterized by)
\begin{equation*}
E(b\ot 1)=E(1\ot S_B(b))
\tussenen
(1\ot c)E=(S_C(c)\ot 1)E
\end{equation*}
for all $b\in B$ and $c\in C$. Then $E$ is called a {\it separability idempotent} in $M(B\ot C)$. 
\edefin

If the maps $S_B$ and $S_C$ have range in $C$ and $B$ respectively, $E$ is called semi-regular and if moreover these maps are anti-isomorphisms, then $E$ is called regular. We will mostly assume that $E$ is regular, although sometimes, we also consider the more general case. 
In the regular case, we have the equality $(S_B\ot S_C)E=\zeta E$ where we use $\zeta$ for the flip map from $B\ot C$ to $C\ot B$, eventually extended to the multiplier algebra $M(B\ot C)$. Consequently, also $(S_CS_B\ot S_BS_C)E=E$.
\snl
The regular case was first studied in the first version of \cite{VD4}, while the more general case is considered in the second version of this paper.
\snl
We denote by $\varphi_B$ and $\varphi_C$ the distinguished linear functionals on $B$ and $C$  satisfying
\begin{equation*}
(\varphi_B\ot \iota)E=1
\tussenen
(\iota\ot\varphi_C)E=1
\end{equation*}
in $M(C)$ and $M(B)$ respectively. Recall that in the regular case, the distinguished functionals are faithful and have KMS automorphisms. The KMS automorphism $\sigma_B$ of $B$, satisfying $\varphi_B(bb')=\varphi_B(b'\sigma_B(b))$ for all $b,b'\in B$ is given by $\sigma_B=(S_CS_B)^{-1}$ while the KMS automorphism $\sigma_C$ of $C$, satisfying $\varphi_C(cc')=\varphi_C(c'\sigma_C(c))$ for all $c,c'\in C$ is given by $\sigma_C=S_BS_C$. Moreover $\varphi_C=\varphi_B\circ S_C$ and $\varphi_B=\varphi_C\circ S_B$.
\snl
We will also need some results here that are not included  the papers \cite{VD4}.
\snl
The following is a {\it Radon-Nikodym type relation} between linear functionals on the base algebras. It is assumed that $E$ is regular for this result.

\prop\label{prop:1.2} 
 Assume that $f$ is a faithful linear functional on $B$. For any other linear functional $g$ on $B$, there is a unique element $y$ in $M(B)$ such that $g(b)=f(by)$ for all $b\in B$. The functional $g$ is faithful if and only if $y$ is invertible in $M(B)$.
\eprop

\bew i) Let $g$ be any linear functional on $B$ and define $y=(\iota\ot (g\circ S_C))E$. This is a well-defined element in $M(B)$. And for all $b\in B$ we get
\begin{align*}
\varphi_B(by)
&=(\varphi_B\ot (g\circ S_C))((b\ot 1)E)\\
&=(\varphi_B\ot (g\circ S_C))((1\ot S_C^{-1}(b))E)\\
&=g(S_CS_C^{-1}(b))=g(b).
\end{align*}
We find that $g(b)=\varphi_B(by)$ for all $b$. 
\snl
If we take $y'=(\iota\ot (g\circ S_B^{-1}))E$, we get another element in $M(B)$ and for all $b$ we have
\begin{align*}
\varphi_B(y'b)
&=(\varphi_B\ot (g\circ S_B^{-1}))(E(b\ot 1))\\
&=(\varphi_B\ot (g\circ S_B^{-1}))(E(1\ot S_B(b))\\
&=g(S_B^{-1}S_B(b))=g(b).
\end{align*}
We now find that $g(b)=\varphi_B(y'b)$ for all $b$. 
\snl
ii) Now assume that $g$ is faithful. We claim that the element $y$ obtained in i) is invertible in $M(B)$. 
\snl
We first show that $By=B$. To prove this, take a linear functional $\omega$ on $B$ with the property that $\omega(by)=0$ for all $b$. By i) there is an element $y'\in M(B)$ so that $\omega(b)=\varphi_B(y'b)$ for all $b\in B$. By assumption $\varphi_B(y'by)=0$ for all $b$. This means that $g(y'b)=0$ for all $b$. As $g$ is assumed to be faithful, this implies that $y'=0$. Hence $\omega=0$. Therefore $By=B$.
\snl
To show that also $yB=B$, let $\omega$ be a linear functional so that $\omega(yb)=0$ for all $b$. Again take $y'\in M(B)$ satisfying $\omega(b)=\varphi_B(y'b)$ for all $b$. Then $\varphi_B(y'yb)=0$ and $\varphi_B(\sigma_B^{-1}(b)y'y)=0$. This means that $g(\sigma_B^{-1}(b)y')=0$ for all $b$. Then $g(by')=0$ for all $b$ and again, by the faithfulness of $g$, that $y'=0$.  So $\omega=0$ and therefore $yB=B$. 

\snl
So we have both $yB=B$ and $By=B$. 
\snl
We will now show that then $y$ is invertible. Define a left multiplier $z$ of $B$ by $zyb=b$ for all $b$. This is well-defined because if $yb=0$, then $b'yb=0$ for all $b'$ and as By=B, also $b'b=0$ for all $b'$ so that $b=0$. Because aslo $yB=B$, the left multiplier $z$ is everywhere defined. Similarly, we have a right multiplier $z'$ satisfying $byz'=b$ for all $b$. Then 
\begin{equation*}
((by)z')(yb')=b(yb')=(by)b'=(by)(z(yb))
\end{equation*}
for all $b,b'$. Using that $By=B$ and $yB=B$, we also get $(bz')b'=b(zb')$ for all $b,b'$. It follows that $z$ is a multiplier in $M(B)$ and $z'=z$.  Now $zyb=b$ and $byz=b$ for all $b$ imply that $yz=1$ and $zy=1$. So $y$ is invertible.
\snl
iii) We are now ready to complete the proof of the proposition. Take two linear functionals $f$ and $g$ on $B$ and assume that $f$ is  faithful. There exist elements $y,y'\in M(B)$ satisfying
\begin{equation*}
f(b)=\varphi_B(by)
\tussenen
g(b)=\varphi_B(by')
\end{equation*}
for all $b$. As $f$ is faithful, the element $y$ has an inverse $z$ in $M(B)$. Then
\begin{equation*}
g(b)=\varphi_B(by')=f(by'z)
\end{equation*}
for all $b$. This proves the first statement. Finally, if also $g$ is faithful, then also $y'$ is invertible and so is $y'z$. This proves the second statement.
\ebew 

We have similar results for linear functionals on the algebra $C$.
\snl
Moreover, we have the following corollary.

\prop
Any faithful linear functional on $B$ has the KMS property.
\eprop

The argument is simple. Take $f$ of the form $b\mapsto \varphi_B(by)$ with $y$ invertible in $M(B)$. Then
\begin{equation*}
f(bb')=\varphi_B(bb'y)=\varphi_B(b'y\sigma_B(b))=f(b'y\sigma_B(b)y^{-1})
\end{equation*}
for all $b,b'\in B$.
\snl
Let us finish this short survey on separability idempotents with a remark. In Theorem 4.13 of \cite{VD4}, a structure theorem is obtained for the base algebras of a regular separability idempotent. The above properties follow also from that result in a more or less straightforward manner.
\nl
\bf Weak multiplier Hopf algebras with integrals \rm
\nl
In \cite{VD-W4} weak multiplier Hopf algebras are introduced as a pair $(A,\Delta)$ where $A$ is a non-degenerate idempotent algebra and $\Delta$ a coproduct on $A$, subject to a number of conditions. The conditions are motivated in a preliminary paper on the subject (\cite{VD-W3}). The source and target maps are introduced already in \cite{VD-W3} but the more complete study is found in \cite{VD-W5}. Integrals on regular weak multiplier Hopf algebras, as well as duality are treated in \cite{VD-W6}. And as we mentioned already in the introduction, this note is meant to illustrate some of the aspects of that paper.
\snl
On the other hand, there is also the recent work on the Larson-Sweedler theorem for weak multiplier Hopf algebras \cite{K-VD}. Roughly speaking, in that paper, it is shown that a weak multiplier bialgebra, as introduced in \cite{B-T-L}, is a regular weak multiplier Hopf algebra if it has (enough) integrals. And since we are interested here in the integrals for our example, and the duality, instead of recalling the general notion as in \cite{VD-W3} of a weak multiplier Hopf algebra, we recall the setting as in \cite{K-VD}. 
\snl
In any case, we start with a pair of a non-degenerate idempotent algebra $A$ and a regular coproduct $\Delta$ on $A$. We assume that it is weakly non-degenerate. More precisely, we have the following.

\defin
A coproduct on $A$ is a homomorphism $\Delta:A\to M(A\ot A)$. We first assume that elements of the form
\begin{align*}
&\Delta(a)(1\ot b) \tussenen (c\ot 1)\Delta(a) \\
&(1\ot b)\Delta(a) \tussenen \Delta(a)(c\ot 1)
\end{align*}
belong to $A\ot A$ for all $a,b,c\in A$. Next we assume that $\Delta$ is coassociative.
\snl
We also assume that $\Delta$ is weakly non-degenerate. This means that there exists an idempotent $E$ in $M(A\ot A)$ with the property that 
\begin{equation*}
\Delta(A)(A\ot A)=E(A\ot A)
\tussenen
(A\ot A)\Delta(A)=(A\ot A)E.
\end{equation*}
\edefin

Because $\Delta$ is weakly non-degenerate, the homomorphisms $\Delta\ot\iota$ and $\iota\ot\Delta$ have unique extension to homomorphisms form $M(A)$ to $M(A\ot A)$ provided we assume 
\begin{equation*}
(\Delta\ot\iota)(1)=E\ot 1
\tussenen
(\iota\ot\Delta)(1)=1\ot E,
\end{equation*}
where we use the same symbols for the extensions of $\Delta\ot\iota$ and $\iota\ot\Delta$. Recall that we use $\iota$ for the identity map. It is easy to show that coassociativity is equivalent with $(\Delta\ot\iota)\Delta=(\iota\ot\Delta)\Delta$. See e.g.\ the short note on coassociativity \cite{VD8}. 

\nl
The next assumption is 
\begin{equation}
(\Delta\ot\iota)(E)=(E\ot 1)(1\ot E)=(1\ot E)(E\ot 1).\label{eqn:Delta-On-E}
\end{equation}
These assumptions are also part of the axioms for a regular weak multiplier Hopf algebra as in Definitions 1.14 and Definition 4.1 of \cite{VD-W4}. And they are also the first assumptions to obtain the Larson-Sweedler theorem in Section 2 of \cite{K-VD}. 
\snl
Now we continue along the lines of \cite{K-VD}. 

\voorw\label{defin:sep-id}
We assume that $E$ is a separability idempotent in the following sense. There exists subalgebras $B$ and $C$ sitting non-degenerately in $M(A)$. This means that $BA=AB=A$ as well as $CA=AC=A$. And it is assumed that $E$ is a separability idempotent in $M(B\ot C)$. 
\evoorw

It is important to observe that the algebras $B$ and $C$ are completely determined by these conditions. This is the reason why we formulate the assumption just by saying that $E$ is a separability idempotent.
\snl
The algebras $B$ and $C$ are automatically non-degenerate and their multiplier algebras $M(B)$ and $M(C)$ sit in $M(A)$ as well. We denote them by $A_s$ and $A_t$ respectively.
\snl
Now we  give the definition of left and right integrals in this setting.

\defin
A left integral on $A$ is a non-zero linear functional $\varphi$ on $A$ with the property that $(\iota\ot\varphi)\Delta(a)\in A_t$ for all $a\in A$. Similarly, a right integral on $A$ is a non-zero linear functional $\psi$ on $A$ with the property that $(\psi\ot\iota)\Delta(a)\in A_s$ for all $a\in A$.
\edefin

We now assume that there is a faithful set of left integrals in the following sense. Suppose that $a\in A$ and that $\varphi(ab)=0$ for all $b$ and all left integrals $\varphi$. Then we must have $a=0$. Similarly, if $\varphi(ca)=0$ for all 
$c\in A$ and all left integrals $\varphi$. We also assume that there is a faithful set of right integrals.
\snl
This takes us to the following definition.

\defin 
Let $(A,\Delta)$ be a pair of a non-degenerate idempotent algebra $A$ and $\Delta$ a weakly non-degenerate coproduct on $A$. Assume that the canonical idempotent $E$ satisfies the equation (\ref{eqn:Delta-On-E}). Assume also that $E$ is a separability idempotent in the sense of Definition \ref{defin:sep-id}. We call $(A,\Delta)$ an {\it algebraic quantum groupoid} if there is a faithful set of left integrals and a faithful set of right integrals.
\edefin

It is shown in Theorem 2.14 of \cite{K-VD} that an algebraic quantum group is nothing else but a regular weak multiplier Hopf algebra with a faithful set of integrals. 
\nl\nl

\section{\hspace{-17pt}. The weak multiplier Hopf algebra of a separability idempotent}\label{s:thewmha} 

In this section, we start with two non-degenerate algebras $B$ and $C$ and a separability idempotent $E$ in $M(B\ot C)$. We recall the construction of the associated weak multiplier Hopf algebra as it is first obtained in Section 3 of\cite{VD-W5}. Also for this paper, there are two different versions on the arXiv. If we don't mention it, internal references are with respect to the second version. We include some of the details, found already in that paper, just for completeness. What is new is the study of integrals on this weak multiplier Hopf algebra.

\nl
\bf The weak multiplier  Hopf algebra of a separability idempotent \rm
\nl
We have the following associated weak multiplier Hopf algebra (cf.\ Proposition 3.2 in \cite{VD-W5}). 

\prop \label{prop:2.1}
Let $A$ be the algebra $C\ot B$. Define $\Delta$ on $A$ by
\begin{equation*}
\Delta(c\ot b)=c\ot E\ot b
\end{equation*}
for all $b\in B$ and $c\in C$. Then $(A,\Delta)$ is a weak multiplier Hopf algebra. 
\vskip 3pt
The counit is given by $\varepsilon(c\ot b)=\varphi_B(S_C(c)b)=\varphi_C(cS_B(b))$. The antipode satisfies $S(c\ot b)=S_B(b)\ot S_C(c)$ and the source and target maps are 
\begin{equation*}
\varepsilon_s(c\ot b)=1\ot S_C(c)b
\tussenen
\varepsilon_t(c\ot b)=cS_B(b)\ot 1
\end{equation*}
for all $b,c$.
\eprop

We see that $\varepsilon_s(A)=1\ot B$ and $\varepsilon_t(A)=C\ot 1$. The multiplier algebras $A_s$ and $A_t$ are $1\ot M(B)$ and $M(C)\ot 1$ respectively. 
\snl
Observe also that $A$ is regular if and only if $E$ is regular (in the sense of Definition 2.3 of \cite{VD4}). See Proposition 3.3 in \cite{VD-W5}.
\snl
In what follows we will look at $A$ as the algebra generated by $B$ and $C$ with the condition that elements of $B$ commute with elements of $C$. We can do this by considering the natural non-degenerate homomorphisms from $B$ and $C$ to $M(A)$ and by identifying $B$ and $C$ with their images. Then $A$ is spanned by elements of the form $cb$ and $cb=bc$ in $A$. This is standard and a common practice. 
\snl
Then we write $\Delta(cb)=(c\ot 1)E(1\ot b)$ and the canonical idempotent of $A$ is $E$, considered as sitting in $M(A\ot A)$. We refer to Section 3 in  \cite{VD-W5}.
\snl
We get the following characterization of the integrals.

\prop\label{prop:2.2}
For any linear functional $g$ on $B$ we have a left integral $\varphi$ on $A$ given by $\varphi(cb)=\varphi_C(c)g(b)$ and  any left integral has this form. Similarly, for any linear functional $f$ on $C$ we have a right integral $\psi$ on $A$ given by $\psi(cb)=f(c)\varphi_B(b)$ and again any right integral is of this form. Moreover, the map $cb\mapsto \varphi_C(c)\varphi_B(b)$ is a faithful left integral that is also a right integral. Recall that $\varphi_B$ and $\varphi_C$ are the distinguished linear functionals on $B$ and $C$ respectively.
\eprop

\bew
i) Let $g$ be any linear functional on $B$ and define $\varphi$ on $A$ by $\varphi(cb)=\varphi_C(c)g(b)$. We claim that $\varphi$ is a left integral on $A$. Indeed, consider $b\in B$ and $c\in C$. Then
\begin{align*}
(\iota\ot\varphi)\Delta(cb)
&=(\iota\ot\varphi)((c\ot 1)E(1\ot b))\\
&=(\iota\ot\varphi_C)((c\ot 1)E)g(b)\\
&=cg(b).
\end{align*}
This belongs to $C$ and as $A_t=M(C)$ in this case, we see that $\varphi$ is left invariant.
\vskip 3pt
ii) Conversely, assume that $\varphi$ is a left integral on $A$. So
\begin{equation*}
(\iota\ot\varphi)((c\ot 1)E(1\ot b))\in M(C)
\end{equation*}
for all $b\in B$ and $c\in C$. When we multiply with $c'\in C$ from the left, we get something in $C$. Then, as $C$ is idempotent, we find that actually 
\begin{equation*}
(\iota\ot\varphi)((c\ot 1)E(1\ot b))\in C
\end{equation*}
for all $b,c$. 
Then we must have $(\iota\ot\varphi)(E(1\ot b))=g(b)1$ for some linear functional $g$ on $B$. Apply $\omega(c\,\cdot\,)$, with $c\in C$ and $\omega$ a linear functional on $C$. Write $c'=(\omega(c\,\cdot\,)\ot\iota)E$. We find that 
\begin{equation*}
\varphi(c'b)=\omega(c)g(b)=\varphi_C(c')g(b).
\end{equation*}
As the left leg of $E$ is all of $C$ we get this for all $c'\in C$.
\vskip 3pt
iii) Similarly, for any linear functional $f$ on $C$ we have a right integral $\psi$ on $A$ defined by $\psi(cb)=f(c)\varphi_B(b)$. And every right integral is of this form.
\vskip 3pt
iv) Finally, it follows that the map $\varphi: cb\mapsto \varphi_C(c)\varphi_B(b)$ defines a linear functional that is left and right invariant. And because the distinguished linear functionals $\varphi_B$ and $\varphi_C$ are faithful, also this linear functional $\varphi$ is faithful on $A$.
\ebew

We can illustrate some of the results we obtained in the case of the existence of a single faithful integral in Section 1 of \cite{VD-W6}.
\snl
It is fairly straightforward to show e.g.\ that the modular automorphism $\sigma_\varphi$ of the integral $\varphi$, defined as $cb\mapsto \varphi_C(c)\varphi_B(b)$, is given by $\sigma_\varphi(cb)=\sigma_C(c)\sigma_B(b)$ where $\sigma_C$ and $\sigma_B$ are the KMS automorphisms of $\varphi_C$ and $\varphi_B$ respectively.
\snl
To illustrate the Propositions 1.8 and 1.9 of \cite{VD-W6}, we will use the  Radon-Nikodym type properties of the distinguished linear functionals $\varphi_B$ and $\varphi_C$. The result, as formulated here, is a special case of the more general result we have included in the preliminary section, see Proposition \ref{prop:1.2}. In fact, the special result below was proven first to obtain the more general result.

\prop For any linear functional $g$ on $B$, there is a unique element $y\in M(B)$ such that $g(b)=\varphi_B(by)$ for all $b\in B$. If $g$ is faithful, then $y$ is invertible. Similarly, for any linear functional $f$ on $C$ there is a unique element $x\in M(C)$ such that $f(c)=\varphi_C(cx)$ for all $c\in C$. Again, if $f$ is faithful, then $x$ is invertible.
\eprop

Now, to illustrate Proposition 1.8 of \cite{VD-W6}, take the faithful integral $\varphi$ defined as before by $cb\mapsto \varphi_C(c)\varphi_B(b)$. Take now any left integral $\varphi_1$. It has the form $\varphi_1(cb)=\varphi_C(c)g(b)$ for some linear functional $g$ on $B$. By the previous proposition, there is an element $y\in M(B)$ so that $g(b)=\varphi_B(by)$ for all $b\in b$. This means that 
\begin{equation*}
\varphi_1(cb)=\varphi_C(c)\varphi_B(by)=\varphi(cby)
\end{equation*}
for all $b\in B$ and $c\in C$. This is the formula we obtained in Proposition 1.8 of \cite{VD-W6}.
\snl
To illustrate Proposition 1.9 of \cite{VD-W6}, consider again the left integral $\varphi$ and right integral $\psi$. The first one is of the form $cb\mapsto \varphi_C(c)g(b)$ for a linear functional $g$ on $B$ and the second one has the form $cb\mapsto f(c)\varphi_B(b)$ for a linear functional $f$ on $C$. Assume that $g$ is faithful so that $\varphi$ is faithful. By the previous proposition, there  elements $x\in M(C)$ and $y\in M(B)$ satisfying $f=\varphi_C(\,\cdot\, x)$ and $g=\varphi_B(\,\cdot\,y)$. Because we assume that $g$ is faithful, the element $y$ has an inverse $y^{-1}$ in $M(B)$. Then
\begin{align*}
\psi(cb)=f(c)\varphi_B(b))=\varphi_C(cx)g(by^{-1})=\varphi(cxby^{-1})=\varphi(cbxy^{-1}).
\end{align*}
We indeed get a modular element $\delta \in M(A)$ satisfying $\psi=\varphi(\,\cdot\, \delta)$ as in Proposition 1.9 of \cite{VD-W6}.

\section{\hspace{-17pt}. Construction of the dual} \label{s:duality} 

As in the previous section, we start with a regular separability idempotent $E$ in $M(B\ot C)$ where $B$ and $C$ are non-degenerate algebras. We consider the algebraic quantum groupoid $(A,\Delta)$ as defined before. Recall that $A$ is the algebra generated by $B$ and $C$ with the condition that elements of $B$ commute with elements of $C$ in $A$. The coproduct $\Delta$ on $A$ is defined by $\Delta(cb)=(c\ot 1) E(1 \ot b)$ for $b\in B$ and $c\in C$.
\snl
Recall that the algebras have local units and in particular, they are idempotent (see Proposition 1.10 in \cite{VD4}). 
\nl
We will now construct the dual $\widehat A$ following the general procedure as described in Section 2 of \cite{VD-W6}.
 In order to do this, we need that $(A,\Delta)$ is a regular weak multiplier Hopf algebra and we know that this is the case if and only if $E$ is a regular separability idempotent. So we assume this in what follows.
\snl
In the regular case, we have a faithful integral $\varphi$ on $A$, defined by $\varphi(cb)=\varphi_C(c)\varphi_B(b)$. Moreover,  as a vector space, $\widehat A$ can be identified with $C\ot B$ via the map 
$v\ot u\mapsto \varphi(\,\cdot\,vu)$. Remark that $\varphi(cbvu)=\varphi_C(cv)\varphi_B(bu)$ for all $b,u\in B$ and $c,v\in C$. So this map is $v\ot u\mapsto \varphi_C(\,\cdot\,v)\ot \varphi_B(\,\cdot\,u)$. 
\snl
In what follows however, we will make another identification. It turns out to give nicer formulas.

\defin\label{defin:pairing}
We define a pairing of the vector spaces $A$ and $B\ot C$ by
\begin{equation*}
\langle cb,u\ot v \rangle = \varphi_B(bS_C(v))\varphi_C(S_B(u)c).
\end{equation*}
\edefin

Observe the difference with the formula for the counit on $A$. Indeed, for the counit we have $\varepsilon(cb)=\varphi_B(S_C(c)b)=\varphi_C(cS_B(b))$, see Proposition \ref{prop:2.1}.
\snl
Further we have
\begin{align*}
 \varphi_B(bS_C(v))\varphi_C(S_B(u)c)
 &=\varphi_B(bS_C(v)\varphi_C(c\sigma_C(S_B(u))) \\
 &=\varphi(cbS_C(v)\sigma_C(S_B(u))).
\end{align*}
We know that $S_C$ and $S_B$ are bijections from $B$ to $C$ and from $C$ to $B$ respectively. And as also $\sigma_C$ is a bijective map from $C$ to itself, we see that $u\ot v\mapsto S_C(v)\ot\sigma_C(S_B(u))$ will be a one-to-one map from $B\ot C$ to $B\ot C$. Therefore 
\begin{equation}
u\ot v\mapsto \varphi(\,\cdot\,S_C(v)\sigma_C(S_B(u)))
\end{equation}
gives a linear isomorphism from the space $B\ot C$ to $\widehat A$.
\snl

In what follows, we will make this identification, but we will use $u\diamond v$ for $u\ot v$ when we consider it as an element in $\widehat A$. We will also systematically use the letters $u,u',\dots$ and $v,v'\dots$ for elements in $B$ and $C$ respectively, when $B\ot C$ is identified with $\widehat A$ via the pairing in Definition \ref{defin:pairing} above.
\nl
\bf The algebra $\widehat A$ and its multiplier algebra $M(\widehat A)$ \rm
\nl
For the product on the dual $\widehat A$, we find the following result.

\prop\label{prop:4.10} 
For all $u,u'\in B$ and $v,v'\in C$ we have
\begin{equation}
(u\diamond v)(u'\diamond v')=\varepsilon(vu')\,\,u\diamond v' \label{eqn:3.5}
\end{equation}
were $\varepsilon$ is the counit on $A$.
\eprop

\bew
We have by the definition of the product in $\widehat A$, for all $b,u,u'$ in $B$  and $c,v,v'$ in $C$,
\begin{align*}
\langle cb,(u\di v)(u'\di v')\rangle
&=\langle (c\ot 1)E(1\ot b),(u\di v)\ot (u'\di v')\rangle \\
&=\langle cE_{(1)},u\di v\rangle \langle E_{(2)}b,u'\di v'\rangle \\
&=\varphi_B(E_{(1)}S_C(v))\varphi_C(S_B(u)c)\varphi_B(bS_C(v'))\varphi_C(S_B(u')E_{(2)})\\
&=\varphi_B(E_{(1)})\varphi_C(S_B(u)c)\varphi_B(bS_C(v'))\varphi_C(S_B(u')E_{(2)}S_BS_C(v))\\
&=\varphi_C(S_B(u)c)\varphi_B(bS_C(v'))\varphi_C(S_B(u')S_BS_C(v))\\
&= \varphi_C(S_B(u')S_BS_C(v))\,\,\langle cb, u\di v'\rangle\\
&=\varphi_B(S_C(v)u')\,\,\langle cb, u\di v'\rangle\\
&=\varepsilon(vu')\,\, \langle cb, u\di v'\rangle.
\end{align*}
We have used the Sweedler type notation $E_{(1)}\ot E_{(2)}$ for $E$ and 
\begin{equation*}
 E_{(1)}S_C(v)\ot E_{(2)}=E_{(1)}\ot E_{(2)}S_B(S_C(v)).
\end{equation*} 
\snl
This proves the formula (\ref{eqn:3.5}) in the formulation of the proposition.
\ebew

By the general theory, we know that the product on $\widehat A$ is associative and non-degenerate. Here, we can verify that the product is associative by a simple calculation. We see also that it is non-degenerate  because the map $(u,v)\mapsto \varepsilon(vu)$ is a non-degenerate pairing of the space $B$ with $C$. This follows from the definition $\varepsilon(vu)=\varphi_B(S_C(v)u)$, the fact that $\varphi_B$ is faithful and that $S_C$ is bijective.
\snl
It is also easy to see that the algebra is idempotent and in fact has local units.
\nl
In what follows, we will write $B\di C$ for the algebra we obtain in Proposition \ref{prop:4.10}. The algebra is an infinite matrix algebra, build with two vector spaces and a non-degenerate pairing. The algebra structure is not dependent on the multiplications in $B$ and $C$. It only depends on the pairing of the underlying vector spaces.
\snl
Because the algebra $B\di C$ is non-degenerate, we can consider the multiplier algebra \newline $M(B\di C)$. We get the following characterization. 

\prop\label{prop:3.3}
A linear map $\gamma:B\to B$ is called adjointable if there is a linear mapping $\gamma^t:C\to C$ satisfying $\varepsilon(v\gamma(u))=\varepsilon(\gamma^t(v)u)$ for all $u\in B$ and $v\in C$. The maps $\gamma$ and $\gamma^t$ determine each other. For any adjointable map $\gamma:B\to B$ with adjoint $\gamma^t$, there is a multiplier $m$ of $B\di C$ given by
\begin{equation}
m(u\di v)=\gamma(u)\di v
\tussenen
(u\di v)m=u\di\gamma^t(v). \label{eqn:multiplier}
\end{equation}
Any multiplier is of this form.
\eprop

\bew
i) For all $u,u'$ in $B$ and $v,v'$ in $C$ we have
\begin{equation*}
(u\di v)(\gamma(u')\di v')
=\varepsilon(v\gamma(u'))\,u\di v'
=\varepsilon(\gamma^t(v)u')\,u\di v'
=(u\di \gamma^t(v))(u'\di v'). 
\end{equation*} 
This proves that a multiplier can be defined as in equation (\ref{eqn:multiplier}).
\snl
ii) Conversely, assume that $m$ is a multiplier of $B\di C$. Then we have 
\begin{equation*}
m((u\di v)(u'\di v'))=(m(u\di v))(u'\di v')
\end{equation*}
for all $u,u'$ and $v,v'$. From the definition of the product, it follows from this that $m(u\di v)$ must be of the form $m(u\di v)=\gamma(u)\di v$ for a linear map $\gamma:B\to B$. Similarly $(u\di v)m$ has the form $u\di \gamma'(v)$ for a linear map $\gamma':C\to C$. The equality $(u\di v)(m(u'\di v'))=((u\di v)m)(u'\di v')$ implies that $\gamma$ and $\gamma'$ are each others adjoints.
\ebew

The embedding of $B\di C$ in $M(B\di C)$ is found by associating the linear map $u'\mapsto \varepsilon(vu')u$ to the element $u\di v$. Its adjoint is $v'\mapsto \varepsilon(v'u)v$. The identity map from $B$ to itself is of course adjointable and the associated multiplier is the identity $1$ in $M(B\di C)$.

\opm
Consider the following argument. We again use the Sweedler type notation $E_{(1)}\ot E_{(2)}$ for $E$. For all $u,v$ we have
\begin{align}
(E_{(1)}\di v)\, \varepsilon(E_{(2)}u)
&=(E_{(1)}\di v)\, \varphi_C(E_{(2)}S_B(u))\label{eqn:4.34}\\
&=(E_{(1)}u\di v)\, \varphi_C(E_{(2)})\label{eqn:4.35}\\
&= u\di v.\nonumber
\end{align}
Formally, we could write  $(E_{(1)}\di E_{(2)})(u\di v)$ for the first expression in this series of equalities and this would express $1$ as the formal expression  $E_{(1)}\di E_{(2)}$. The above calculation shows that this is well-defined as a map from $B$ to itself because the factor $E_{(2)}$ is covered by $S_B(u)$ in (\ref{eqn:4.34}) while the factor $E_{(1)}$ is covered by $u$ in (\ref{eqn:4.35}).
\snl
This result seems interesting and certainly instructive, but we will not really use it explicitly.
\eopm
 
Remark again that the algebra structure of $B\di C$ is only dependent on the vector spaces $B$ and $C$, together with the pairing $u\ot v\mapsto \varepsilon(uv)$. This is not so for the coproduct we consider in the next item.
\nl
 \bf The coproduct on $B\di C$\rm
 \nl
 By definition, the coproduct on $\widehat A$ is dual to the product on $A$. Because $A$ is the tensor product of the algebras $C$ and $B$, it is expected that the coproduct on the dual is also a tensor product of coproducts on the factors. In this case, it would mean that the coproduct on $B\di C$ has the form 
 \begin{equation*}
\Delta(u\di v)=\sum_{(u),(v)} (u_{(1)}\di v_{(1)})\ot (u_{(2)}\di v_{(2)})
\end{equation*}
where we have the Sweedler notation for coproducts $\Delta_B$ on $B$ and $\Delta_C$ on $C$.
\snl
This is indeed the case and it is made precise with the following results.
\snl
Recall the notations $F_1=(\iota\ot S_C)E$ and $F_2=(S_B\ot\iota)E$ and the  formulas
\begin{equation*}
(u\ot 1)F_1=F_1(1\ot u)
\tussenen
(v\ot 1)F_2=F_2(1\ot v)
\end{equation*}
for $u\in B$ and $v\in C$. These formulas are easy to prove. See e.g.\ Section 2 of \cite{VD4}.

\prop\label{prop:3.5} 
Define a pairing of $C$ with $B$ by the formula $\langle c,u \rangle_1=\varphi_C(S_B(u)c)$. For this pairing we have
\begin{equation*}
\langle cc',u \rangle_1 = \langle c\ot c',\Delta_B(u) \rangle_1
\end{equation*}
for $c,c'\in C$ and $u\in B$, where $\Delta_B(u)=F_1(1\ot u)$.
\eprop

\bew
Take $u\in B$ and $c,c'\in C$. Then
\begin{align*}
\langle cc',u \rangle_1
&=\varphi_C(S_B(u)cc')\\
&=(\varphi_B\ot\varphi_C)((1\ot S_B(u)c)E(1\ot c')) \\
&=(\varphi_B\ot\varphi_C)((S_C(c)\ot S_B(u))E(1\ot c')) \\
&=(\varphi_B\ot\varphi_C)((S_C(c)\ot S_B(u))((S_CS_B\ot S_BS_C)E)(1\ot c')) \\
&=(\varphi_C\ot\varphi_C)((1\ot S_B(u))((S_B\ot S_BS_C)E)(c\ot c')) \\
&=(\varphi_C\ot\varphi_C)((S_B\ot S_B)(((\iota\ot S_C)E)(1\ot u))(c\ot c')) \\
&=\langle F_1(1\ot u),c\ot c'\rangle _1
\end{align*}
\vskip -0.8cm
\ebew

\prop\label{prop:3.6}
Define a pairing of $B$ with $C$ by the formula $\langle b,v \rangle_2= \varphi_B(bS_C(v))$. For this pairing we have 
\begin{equation}
\langle bb',v\rangle_2 =\langle b\ot b',\Delta_C(v)\rangle_2
\end{equation}
where $\Delta_C(v)=(v\ot 1)F_2$ for $b,b'\in B$ and $v\in C$.
\eprop
 
\bew
 Take $b,b'\in B$ and $v\in C$. Then
 \begin{align*}
\langle bb',v\rangle_2
&=\varphi_B(bb'S_C(v)) \\
&=(\varphi_B\ot\varphi_C)((b\ot 1)E(b'S_C(v)\ot 1))\\
&=(\varphi_B\ot\varphi_C)((b\ot 1)E(S_C(v)\ot S_B(b')))\\
&=(\varphi_B\ot\varphi_C)((b\ot 1)((S_CS_B\ot S_BS_C)E)(S_C(v)\ot S_B(b')))\\
&=(\varphi_B\ot\varphi_B)((b\ot b')((S_CS_B\ot S_C)E)(S_C(v)\ot 1))\\
&=(\varphi_B\ot\varphi_B)((b\ot b')(S_C\ot S_C)((v\ot 1)(S_B\ot\iota)E))\\
&=(\varphi_B\ot\varphi_B)((b\ot b')(S_C\ot S_C)((v\ot 1)F_2))\\
&=\langle b\ot b',(v\ot 1)F_2\rangle_2.
\end{align*}
 \vskip -0.8cm
\ebew

Now we get from these two properties the formula for the coproduct on $B\di C$.

\prop\label{prop:3.7}
 The coproduct on $B\di C$, defined by the formula 
\begin{equation*}
\langle aa',u\di v\rangle=\langle a\ot a',\Delta(u\di v)\rangle
\end{equation*}
for $a,a'\in A$ and $u\in B$ and $v\in C$ is
\begin{equation*}
\Delta(u\di v)=\sum_{(u),(v)} (u_{(1)}\di v_{(1)})\ot (u_{(2)}\di v_{(2)})
\end{equation*}
where we use the Sweedler notations
\begin{align*}
\sum_{(u)} u_{(1)}\ot u_{(2)}&=\Delta_B(u)=F_1(1\ot u)\\
\sum_{(v)} v_{(1)}\ot v_{(2)}&=\Delta_C(v)=(v\ot 1)F_2.
\end{align*}
\eprop 

We need to make a couple of remarks about this formula.

\opm
i) It is rather remarkable that the pairings $(c,u)\mapsto \langle c , u\rangle_1$ and $(b,v)\mapsto \langle b, v\rangle_2$, that we defined in the Propositions \ref{prop:3.5} and \ref{prop:3.6}, yield coproducts $\Delta_B:B\to B\ot B$ and $\Delta_C:C\to C\ot C$. 
\vskip 2pt
ii) Moreover, these coproducts are known. See e.g.\ Proposition 2.9  in  \cite{VD4}. 
\vskip 2pt
iii) It is also obvious that these coproducts are homomorphisms on the algebras $B$ and $C$ respectively. As a consequence, the tensor coproduct, as defined in Proposition \ref{prop:3.7} will be a homomorphism on the tensor product algebra $B\ot C$. However, this is different from the product in $B\di C$ and, from the general theory, we know that this tensor coproduct is also a homomorphism on $B\di C$.  
\eopm

We will  verify this last remark in the following proposition.

\prop
For all $u,u'\in B$ and $v,v'\in C$ we have
\begin{equation}
\Delta(u\di v)\Delta(u'\di v')=\Delta((u\di v)(u'\di v')) .\label{eqn:3.16}
\end{equation}
\eprop

\bew
We apply the flip map $\zeta$ on the middle two factors in each of the expressions in (\ref{eqn:3.16}).
\snl
For the left hand side we find
\begin{align*}
((u_{(1)}\di  v_{(1)})\ot & (u_{(2)}\di v_{(2)}))((u'_{(1)}\di  v'_{(1)})\ot  (u'_{(2)}\di v'_{(2)}))\\
&=\varepsilon(v_{(1)}u'_{(1)})\varepsilon(v_{(2)}u'_{(2)})(u_{(1)}\di  v'_{(1)})\ot  (u_{(2)}\di v'_{(2)}).
\end{align*}
(We have omitted the summation sign, as we do also in the next series of formulas.)
\snl
We claim that $\varepsilon(v_{(1)}u'_{(1)})\varepsilon(v_{(2)}u'_{(2)})=\varepsilon(vu')$. Then the above expression is precisely $\varepsilon(vu')\zeta_{23}\Delta(u\di v')$ and this will give the equality we want to prove.
\snl
Now, from the definitions of the coproducts $\Delta_B$ and $\Delta_C$ we find
\begin{align*}
\varepsilon(v_{(1)}u'_{(1)})\varepsilon(v_{(2)}u'_{(2)})
&=\varepsilon(vS_B(E_{(1)})E'_{(1)})\varepsilon(E_{(2)}S_C(E'_{(2)})u')\\
&=\varphi_C(vS_B(E_{(1)})S_B(E'_{(1)}))\varphi_B(S_C(E_{(2)})S_C(E'_{(2)})u')\\
&=\varphi_C(vS_B(E'_{(1)}E_{(1)}))\varphi_B(S_C(E'_{(2)}E_{(2)})u')\\
&=\varphi_C(vS_B(E_{(1)}))\varphi_B(S_C(E_{(2)})u')\\
&=\varphi_C(vS_B(u'E_{(1)}))\varphi_B(S_C(E_{(2)}))\\
&=\varphi_C(vS_B(u'))=\varepsilon(vu').
\end{align*}
We have used that $E$ is an idempotent and that $F_1(1\ot u)=(u\ot 1)F_1$.
\snl
This proves the claim and it completes the argument.
\ebew

Before we continue here, observe that the coproduct $\Delta$ maps $\widehat A$ into the tensor product $\widehat A\ot \widehat A$. Still, we do not have an ordinary weak Hopf algebra because the algebra $\widehat A$ does not have an identity (in general). We have a similar situation in the case of a groupoid algebra of a groupoid with infinitely many units. Still, this situation seems to be rather special.
\snl

The following part is devoted to the verification of the various properties that we get from the general result, namely that $B\di C$ with this coproduct, is a regular weak multiplier Hopf algebra. 
\snl
We first look at the counit and the antipode.
\nl
\bf The counit on the algebra $B\di C$ \rm
\nl
Because the coproduct $\Delta$ on $B\di C$ is the tensor product of the coproduct $\Delta_B$ on $B$ with $\Delta_C$ on $C$, it is obvious that the counit $\widehat\varepsilon$ is given by $\widehat\varepsilon(u\di v)=\varepsilon_B(u)\varepsilon_C(v)$ where $\varepsilon_B$ and $\varepsilon_C$ are the counits on $(B,\Delta_B)$ and $(C,\Delta_C)$ respectively. And because we have 
\begin{equation*}
\Delta_B(u)=F_1(1\ot u)
\tussenen
\Delta_C(v)=(v\ot 1)F_2,
\end{equation*}
where $F_1=(\iota\ot S_C)E$ and $F_2=(S_B\ot\iota)E$, we have $\varepsilon_B=\varphi_B$ and $\varepsilon_C=\varphi_C$.
\snl
We now verify that this is in accordance with the general theory. In Proposition 2.10 of \cite{VD-W6} we find that $\widehat\varepsilon$ on $\widehat A$ is given by $\widehat\varepsilon(\omega)=\omega(1)$ for $\omega\in \widehat A$. Because we have the correspondence $u\di v$ with the functional $cb\mapsto \varphi_B(bS_C(v))\varphi_C(S_B(u)c)$, we see that the counit will satisfy $\widehat\varepsilon(u\di v)=\varphi_B(S_C(v))\varphi_C(S_B(u))=\varphi_B(u)\varphi_C(v)$. We indeed find the same result.
\snl

Because the coproduct maps into the algebraic tensor product, the existence of a counit implies that the coproduct is full.

\nl
\bf The antipode \rm
\nl

It is easy to obtain the antipode on the dual. We get it in the following proposition.

\prop\label{prop:3.10}
The antipode on the dual is given by the formula
\begin{equation*}
S(u\di v)=S_B^{-1}(v)\di S_C^{-1}(u)
\end{equation*}
for $u\in B$ and $v\in C$.
\eprop
\bew
Take $u\in B$ and $v\in C$. Then we find
\begin{align*}
\langle cb, S(u\di v)\rangle
&=\langle S(cb),u\di v \rangle=\langle S_B(b)S_C(c),u\di v\rangle \\
&=\varphi_B(S_C(c)S_C(v))\varphi_C(S_B(u)S_B(b))\\
&=\varphi_C(vc)\varphi_B(bu)
\end{align*}
and we see that $S(u\di v)=S_B^{-1}(v)\di S_C^{-1}(u)$.
\ebew

Let us now verify a couple of properties of this antipode.

\prop
The antipode is an anti-isomorphism.
\eprop
\bew
Take $u,u'\in B$ and $v,v'\in C$. Then we have
\begin{equation*}
S((u\di v)(u'\di v'))=\varepsilon(vu')S(u\di v')=\varepsilon(vu')S_B^{-1}(v')\di S_C^{-1}(u).
\end{equation*}
On the other hand
\begin{align*}
S(u'\di v')S(u\di v)
&=(S_B^{-1}(v')\di S_C^{-1}(u'))(S_B^{-1}(v)\di S_C^{-1}(u))\\
&=\varepsilon (S_B^{-1}(v) S_C^{-1}(u'))(S_B^{-1}(v')\di S_C^{-1}(u))
\end{align*}
and because 
\begin{equation*}
\varepsilon (S_B^{-1}(v) S_C^{-1}(u'))=\varphi_B(u'S_B^{-1}(v))=\varphi_C(vS_B(u'))= \varepsilon(vu')
\end{equation*}
we see that $S((u\di v)(u'\di v'))=S(u'\di v')S(u\di v)$. This proves the result. 
\ebew

\prop
The antipode is an anti-coisomorphism.
\eprop
\bew
Take $u\in B$ and $v\in C$. We have 
\begin{equation*}
\Delta(S(u\di v))=\Delta(S_B^{-1}(v)\di S_C^{-1}(u))
=\zeta_{23}(\Delta_B(S_B^{-1}(v))\di \Delta_C(S_C^{-1}(u)))
\end{equation*}
where $\zeta_{23}$ flips the middle two factors.
Using the Sweedler type notation for two copies of $E$, we have
\begin{align*}
\Delta_B(S_B^{-1}(v))&=S_B^{-1}(v)E_{(1)}\ot S_C(E_{(2)}) \\
\Delta_C(S_C^{-1}(u))&=S_B(E'_{(1)})\ot E'_{(2)}S_C^{-1}(u).
\end{align*}
So
\begin{equation*}
\Delta(S(u\di v))=(S_B^{-1}(v)E_{(1)}\di S_B(E'_{(1)})) \ot (S_C(E_{(2)})\di E'_{(2)}S_C^{-1}(u)).
\end{equation*}
On the other hand, 
\begin{equation*}
\Delta(u\di v)=(uE'_{(1)}\di S_B(E_{(1)})) \ot (S_C(E'_{(2)}) \di E_{(2)}v)
\end{equation*}
and therefore 
\begin{equation*}
(S\ot S)\Delta(u\di v)=(E_{(1)}\di S_C^{-1}(uE'_{(1)}))\ot (S_B^{-1}(E_{(2)}v)\di E'_{(2)}).
\end{equation*}
Becasue  $S_B$ and $S_C$ are anti-isomorphisms and  $(S_B\ot S_C)(E)=\zeta E$ we will indeed find
 $\Delta(S(u\di v))=\zeta(S\ot S)\Delta(u\di v)$.  
This completes the argument.
\ebew

\nl
\bf The canonical idempotent for the dual \rm
\nl
Next we look for the element $\widehat E$ in the multiplier algebra $M(\widehat A\ot\widehat A)$ and discuss some of its properties. 
\snl
First remark that we do not expect $\widehat E\in \widehat A\ot\widehat A$, even in this case where the coproduct maps into the tensor product. The reason is that $1$ is not an element of the algebra $\widehat A$.
\snl
In the more easy, finite-dimensional case, we know that $\widehat E=\Delta(1)$ and that $1=E_{(1)}\di E_{(2)}$. It follows that
\begin{align}
\widehat E
&=\zeta_{23}(\Delta_B(E_{(1)}\di \Delta_C(E_{(2)})\nonumber\\
&=\zeta_{23}((E_{(1)}\ot 1)F_1)\di (F_2(1\ot E_{(2)})))\nonumber\\
&=\zeta_{23}(F_1(1\ot E_{(1)}))\di ((E_{(2)}\ot 1)F_2)\label{eqn:3.17}
\end{align}
where as before $F_1=(\iota\ot S_C)E$ and $F_2=(S_B\ot\iota)E$.
\nl

For a precise and correct treatment of $\widehat E$, we will use the formula $\langle a\ot a',\widehat E\rangle=\varepsilon(aa')$, proven in Proposition 2.13 of \cite{VD-W6}. 
\snl
This formula uses the extension of the pairing from $A$ with $\widehat A$ to $A$ with $M(\widehat A)$, in fact more generaly, from the the pairing on $A\ot A$ with $\widehat A\ot \widehat A$ to the pairing of $A\ot A$ with the multiplier algebra $M(\widehat A\ot\widehat A)$. 
\snl
Before we proceed, let us verify that the expression we find for $\widehat E$ in formula (\ref{eqn:3.17}) satisfies this equality. 
Take $b,b'\in B$ and $c,c'\in C$. Then we find
\begin{align*}
\langle cb\ot c'b',&\zeta_{23}(F_1(1\ot E_{(1)}))\di ((E_{(2)}\ot 1)F_2)\rangle \\
&=\langle cb ,E'_{(1)}\di E_{(2)}S_B(E''_{(1)}) \rangle \langle c'b',S_C(E'_{(2)}) E_{(1)}\di E''_{(2)} \rangle \\
&=\varphi_B(bS_C(E_{(2)}S_B(E''_{(1)})))\varphi_C(S_B(E'_{(1)})c) \\
&\qquad\qquad \varphi_B(b'S_C(E''_{(2)}))\varphi_C(S_B(S_C(E'_{(2)})E_{(1)})c') \\
&=\varphi_B(bS_C(E_{(2)}S_C^{-1}(b')))\varphi_C(S_B(S_B^{-1}(c)E_{(1)})c') \\
&=\varphi_B(bb'S_C(E_{(2)}))\varphi_C(S_B(E_{(1)})cc') \\
&=\varphi_B(bb'S_B^{-1}(cc'))=\varphi_C(cc'S_B(bb'))=\varepsilon(cbc'b')
\end{align*}
\snl

Let us now move to a more precise approach, valid not only in the finite-dimensional situation, but in general. 
\snl
In Proposition \ref{prop:3.3} we have seen that an element of the multiplier algebra $M(\widehat A)$ is given by an adjointable linear map from $B$ to itself. Similarly, an element of the multiplier algebra $M(\widehat A\ot\widehat A)$ is given by an adjointable linear map from $B\ot B$ to itself. We will need explicit formulas here in order to find the dual canonical idempotent $\widehat E$ from the formula $\langle a\ot a',\widehat E\rangle=\varepsilon(aa')$. This is what we will do next.
\snl
We claim that the following result is true.

\prop\label{prop:3.10b}
The maps $\gamma: u\ot u'\mapsto (uu'\ot 1) F_1$ from $B\ot B$ to itself and $\gamma':v\ot v'\mapsto F_2(1\ot vv')$ from $C\ot C$ to itself are adjoint to each other (with respect to the tensor product pairing). They are projection maps and they represent the canonical idempotent $\widehat E$ in $M(\widehat A\ot\widehat A)$ (in the sense of Proposition \ref{prop:3.3}).
\eprop

We will prove these properties in a couple of steps. First we have the following lemma.

\lem
For all $u,u'\in B$ and $v,v'\in C$ we have
\begin{equation*}
(\varepsilon\ot\varepsilon)((\gamma(u\ot u'))(v\ot v'))=(\varepsilon\ot\varepsilon)((u\ot u')\gamma'(v\ot v')).
\end{equation*}
\elem

\bew
For the left hand side, using the Sweedler type notation for $E$, we get
\begin{align*}
(\varepsilon\ot\varepsilon)((\gamma(u\ot u'))(v\ot v'))
&=\varepsilon(uu'E_{(1)}v)\varepsilon(S_C(E_{(2)}v')\\
&=\varphi_C(vS_B(uu'E_{(1)}))\varphi_B(S_C(v')S_C(E_{(2)}))\\
&=\varphi_C(vS_B(E_{(1)})S_B(uu'))\varphi_C(E_{(2)}v')\\
&=\varphi_C(vv'S_B(uu'))=\varepsilon(vv'uu').
\end{align*}
Similarly, for the right hand side we find
\begin{align*}
(\varepsilon\ot\varepsilon)((u\ot u')\gamma'(v\ot v'))
&=\varepsilon(uS_B(E_{(1)}))\varepsilon(u'E_{(2)}vv')\\
&=\varphi_C(S_B(E_{(1)})S_B(u))\varphi_B(S_C(E_{(2)}vv')u')\\
&=\varphi_B(uE_{(1)})\varphi_B(S_C(vv')S_C(E_{(2)})u')\\
&=\varphi_B(S_C(vv')uu')=\varepsilon(vv'uu').
\end{align*}•
This proves the result of lemma.
\ebew

That the maps are idempotent follows simply from the fact that $m_BF_1=E_{(1)}S_C(E_{(2)})=1$ as well as $m_CF_2=S_B(E_{(1)})E_{(2)}$. Here we use $m_B$ and $m_C$ for the multiplication maps on $B\ot B$ and $C\ot C$ respectively.
\snl
To complete the proof of Proposition \ref{prop:3.10b} above, we also need to find a formula for the extended pairing. This is obtained in the following lemma.

\lem\label{lem:3.12}
If $\gamma$ is a linear map from $B$ to itself and if $m$ is the associated left multiplier of $B\di C$, then 
\begin{equation*}
\langle cS_CS_B(b),m\rangle=\varphi_C(S_B(\gamma(b))c)
\end{equation*}
for all $b\in B$ and $c\in C$.
\elem
\bew
We have 
\begin{equation*}
\langle cb,(u'\di v')(u\di v)\rangle=\varepsilon(v'u)\varphi_B(bS_C(v))\varphi_C(S_B(u'))c)
\end{equation*}
and because
\begin{align*}
\langle cS_CS_B(u),u'\di v'\rangle
&=\varphi_B(S_CS_B(u)S_C(v'))\varphi_C(S_B(u')c)\\
&=\varphi_B(S_C(v')u)\varphi_C(S_B(u')c)\\
&=\varepsilon(v'u)\varphi_C(S_B(u')c),
\end{align*}
it follows that $(u\di v)\tr (cb)= \varphi_B(bS_C(v))cS_C(S_B(u))$.
\snl
We now pair this formula with $m$. We find
\begin{equation*}
\langle cb, m(u\di v)\rangle= \varphi_B(bS_C(v))\langle cS_C(S_B(u)),m\rangle
\end{equation*}
and the left hand side is
\begin{equation*}
\langle cb, \gamma(u)\di v\rangle=\varphi_B(bS_C(v))\varphi_C(S_B(\gamma(u))c).
\end{equation*}
So 
\begin{equation*}
\varphi_B(bS_C(v))\langle cS_C(S_B(u)),m\rangle=\varphi_B(bS_C(v))\varphi_C(S_B(\gamma(u))c)
\end{equation*}
and  $\langle cS_C(S_B(u)),m\rangle=\varphi_C(S_B(\gamma(u))c)$. This completes the proof.
\ebew

We are now ready to complete the proof of  Proposition  \ref{prop:3.10b}. 
\snl
\bew
We want to show that $\widehat E$ is the multiplier of $\widehat A\ot \widehat A$ associated with the linear map $\gamma: B\ot B\to B\ot B$ given by $\gamma(u\ot u')=(uu'\ot 1)F_1$. Take $b_1,b_2\in B$ and $c_1,c_2\in C$. Using the result of Lemma \ref{lem:3.12} twice we find (with this proposed expression for $\widehat E$)
\begin{align*}
\langle c_1b_1\ot &c_2b_2,\widehat E\rangle\\
&=(\varphi_C\ot\varphi_C)((S_B\ot S_B)(\gamma((S_B^{-1}S_C^{-1}\ot S_B^{-1}S_C^{-1})(b_1\ot b_2)))(c_1\ot c_2))\\
&=(\varphi_C\ot\varphi_C)((S_B\ot S_B)((S_B^{-1}S_C^{-1}(b_1b_2)\ot 1)F_1)(c_1\ot c_2))\\
&=(\varphi_C\ot\varphi_C)(((S_B\ot S_BS_C)E)((S_C^{-1}(b_1b_2)c_1\ot c_2)))\\
&=(\varphi_C\ot\varphi_C)(((S_C^{-1}\ot\iota)E)((S_C^{-1}(b_1b_2)c_1\ot c_2))\\
&=(\varphi_B\ot\varphi_C)(((S_C(c_1)b_1b_2)\ot 1)E)(1\ot c_2))\\
&=\varphi_C((S_C^{-1}(b_1b_2)c_1 c_2))\\
&=\varphi_B(S_C(c_1c_2)b_1b_2)=\varepsilon(c_1c_2b_1b_2)=\varepsilon(c_1b_1c_2b_2).
\end{align*}
This indeed completes the proof of Proposition \ref{prop:3.10b}.
\ebew

We now verify the following result.

\prop
We have 
\begin{equation*}
\widehat E\Delta(u\di v)=\Delta(u\di v)
\tussenen
\Delta(u\di v)=\Delta(u\di v)\widehat E
\end{equation*}
for all $u\in B$ and $v\in C$.
\eprop

\bew
Take $u\in B$ and $v\in C$. 
\snl
i) Using the Sweedler notations for the coproducts $\Delta_B$ and $\Delta_C$ as before, we get
\begin{align*}
\widehat E\Delta(u\di v)
&=\widehat E((u_{(1)}\di v_{(1)})\ot (u_{(2)}\di v_{(2)}))\\
&=(u_{(1)}u_{(2)}E_{(1)}\di v_{(1)})\ot (S(E_{(2)})\di v_{(2)}).
\end{align*}
Now we have $u_{(1)}u_{(2)}=uE_{(1)}S(E_{(2)})=u$
and therefore
\begin{align*}
\widehat E\Delta(u\di v)
&=(uE_{(1)}\di v_{(1)})\ot (S(E_{(2)})\di v_{(2)})\\
&=(u_{(1)}\di v_{(1)})\ot (u_{(2)}\di v_{(2)})\\
&=\Delta(u\di v).
\end{align*}
ii) Similarly we have
\begin{align*}
\Delta(u\di v)\widehat E
&=((u_{(1)}\di v_{(1)})\ot (u_{(2)}\di v_{(2)}))\widehat E\\
&=(u_{(1)}\di S(E_{(1)})\ot (u_{(2)}\di (E_{(2)}v_{(1)}v_{(2)})\\
&=(u_{(1)}\di S(E_{(1)})\ot (u_{(2)}\di (E_{(2)}v)\\
&=(u_{(1)}\di v_{(1)})\ot (u_{(2)}\di v_{(2)})\\
&=\Delta(u\di v).
\end{align*}
Here we have used that $v_{(1)}v_{(2)}=S(E_{(1)})E_{(2)}v=v$.
\ebew
Again, we have omitted the summation signs here.
\nl
\bf More properties of $\widehat E$ \rm
\nl
In the next proposition, we prove some useful formulas. 

\prop\label{prop:3.17}
For all $u\in B$ and $v\in C$ we have
\begin{align}
\widehat E(1\ot (u\di v))&=(E'_{(1)}uE_{(1)}\di E'_{(2)})\ot (S_C(E_{(2)})\di v) \label{eqn:3.18} \\
(\widehat E(u\di v)\ot 1)&=(uE'_{(1)}E_{(1)}\di v) \ot (S_C(E_{(2)})\di E'_{(2)})  \label{eqn:3.19}\\
((u\di v)\ot 1)\widehat E&=(u\di S_B(E_{(1)})) \ot  (E'_{(1)}\di E_{(2)}vE'_{(2)}) \label{eqn:3.20}\\
(1\ot (u\di v))\widehat E&=(E'_{(1)}\di S_B(E_{(1)}))\ot (u\di E_{(2)}E'_{(2)}v). \label{eqn:3.21}
\end{align}
We are using the Sweedler type notation for two copies of $E$.
\eprop

\bew
From the general theory, we know that the expressions on the left are all sitting in $\widehat A\ot\widehat A$, and not just in the multiplier algebra $M(\widehat A\ot\widehat A)$. On the other hand, in the expressions on the right hand side of the four equalities, we have the necessary coverings of the legs of the elements $E$. In the Equations (\ref{eqn:3.18}) and (\ref{eqn:3.19}), the element $u$ covers $E'_{(1)}$ and consequently, also $E_{(1)}$ is covered. In the the Equations (\ref{eqn:3.20}) and (\ref{eqn:3.21}), it is the element $v$ that takes care of the coverings. Hence also the expressions on the right of these equations, are well-defined elements in $\widehat A\ot\widehat A$.
\snl
i) To prove (\ref{eqn:3.18}) we multiply the expression with $u'\di v'$ in the first factor from the right. We find
\begin{align*}
(E'_{(1)}uE_{(1)}\di E'_{(2)})&(u'\di v') \ot (S_C(E_{(2)})\di v) \\
&=\varepsilon(u'E'_{(2)})(E'_{(1)}uE_{(1)}\di v')\ot (S_C(E_{(2)})\di v)\\
&=\varphi_C(E'_{(2)}S_B(u'))(E'_{(1)}uE_{(1)}\di v')\ot (S_C(E_{(2)})\di v)\\
&=(u'uE_{(1)}\di v')\ot (S_C(E_{(2)})\di v)\\
&=\widehat E((u\di v)\ot (u'\di v')).
\end{align*}
For the last equality, we have used the characterization of $\widehat E$ as given in Proposition \ref{prop:3.10b}.
This proves the first formula.
\snl
ii) For (\ref{eqn:3.19}) we multiply again with $u'\di v'$ from the right, but now in the second factor. We get
\begin{align*}
(uE'_{(1)}E_{(1)}\di v) &\ot (S_C(E_{(2)})\di E'_{(2)})(u'\di v') \\
&= \varepsilon(u'E'_{(2)})(uE'_{(1)}E_{(1)}\di v) \ot (S_C(E_{(2)})\di v') \\
&=\varphi_C(E'_{(2)}S_B(u'))(uE'_{(1)}E_{(1)}\di v) \ot (S_C(E_{(2)})\di v') \\
&=(uu'E_{(1)}\di v) \ot (S_C(E_{(2)})\di v') \\
&=\widehat E((u'\di v')\ot (u\di v)).
\end{align*}
This proves the second formula.
\snl
iii) To prove (\ref{eqn:3.20}) we multiply with $u'\di v'$ from the left, in the second factor. We get
\begin{align*}
(u\di S_B(E_{(1)}))& \ot (u'\di v') (E'_{(1)}\di E_{(2)}vE'_{(2)})\\
&=\varepsilon(v'E'_{(1)})(u\di S_B(E_{(1)}))\ot (u'\di E_{(2)}vE'_{(2)})\\
&=\varphi_B(S_C(v')E'_{(1)})(u\di S_B(E_{(1)}))\ot (u'\di E_{(2)}vE'_{(2)})\\
&=(u\di S_B(E_{(1)}))\ot (u'\di E_{(2)}vv')\\
&=((u\di v)\ot (u'\di v'))\widehat E
\end{align*}
where for the last equality, we again refer to Proposition \ref{prop:3.10b}. This proves the third formula.
\snl
iv) Finally, to prove (\ref{eqn:3.21}), multiply the expression on the right with $u'\di v'$ from the left in the first factor. Then
\begin{align*}
(u'\ot v')(E'_{(1)}\di &S_B(E_{(1)}))\ot (u\di E_{(2)}E'_{(2)}v) \\
&=\varepsilon(v'E'_{(1)})(u'\di  S_B(E_{(1)}))\ot (u\di E_{(2)}E'_{(2)}v) \\
&=\varphi_B(S_C(v')E'_{(1)})(u'\di S_B(E_{(1)}))\ot (u\di E_{(2)}E'_{(2)}v) \\
&=(u'\di S_B(E_{(1)}))\ot (u\di E_{(2)}v'v) \\
&=((u'\di v')\ot (u\di v))\widehat E.
\end{align*} 
This completes the proof.
\ebew

From these formulas, we can easily find the legs of $\widehat E$. Indeed, consider e.g.\ Equation  (\ref{eqn:3.18}) and apply a linear functional on the second leg of this equation. This results in elements of the form $E'_{(1)}u_1\di E'_{(2)}$ where $u_1$ can be any element of $B$. Similarly, if we apply a linear functional on the first leg of Equation (\ref{eqn:3.20}), we end up with elements of the form $E'_{(1)}\di v_1E'_{(2)}$ where now $v_1$ can be any element of $C$.
\snl
This result will also be a consequence of the formulas for the source and target maps we consider next.
\nl
\bf The source and target maps \rm
\nl
Using the expression for the antipode, we find the following results for the source and target maps.

\prop
For all $u\in B$ and $v\in C$ we have
\begin{equation*}
\varepsilon_s(u\di v)=\varphi_B(u)(E_{(1)}\di E_{(2)}v)
\tussenen
\varepsilon_t(u\di v)=\varphi_C(v)(uE_{(1)}\di E_{(2)}).
\end{equation*}
\eprop
\bew 
We have 
\begin{align*}
\varepsilon_s(u\di v)
&=S(uE_{(1)}\di S_B(E'_{(1)}))(S_C(E_{(2)})\di E'_{(2)}v)\\
&=(E'_{(1)}\di S_C^{-1}(uE_{(1)}))(S_C(E_{(2)})\di E'_{(2)}v)\\
&= \varepsilon (S_C^{-1}(uE_{(1)})S_C(E_{(2)}))(E'_{(1)}\di E'_{(2)}v)\\
&= \varphi_B (uE_{(1)}S_C(E_{(2)}))(E'_{(1)}\di E'_{(2)}v)\\
&= \varphi_B (u)(E'_{(1)}\di E'_{(2)}v).
\end{align*}
This proves the first equality. For the second one we have
\begin{align*}
\varepsilon_t(u\di v)
&=(uE_{(1)}\di S_B(E'_{(1)}))S(S_C(E_{(2)})\di E'_{(2)}v)\\
&=(uE_{(1)}\di S_B(E'_{(1)}))(S_B^{-1}( E'_{(2)}v)\di E_{(2)})\\
&=\varepsilon(S_B(E'_{(1)})S_B^{-1}( E'_{(2)}v))(uE_{(1)}\di  E_{(2)})\\
&=\varphi_C(S_B(E'_{(1)})E'_{(2)}v)(uE_{(1)})\di  E_{(2)})\\
&=\varphi_C(v)(uE_{(1)})\di  E_{(2)}).
\end{align*}
\ebew

We can also use the formulas of Proposition \ref{prop:3.17} to obtain these results. Indeed, if we apply the counit on the first leg of Equation (\ref{eqn:3.19}) we find
\begin{align*}
\varepsilon_t(u\di v)
&=(\widehat\varepsilon\ot\iota)(\widehat E(u\di v)\ot 1))\\
&=\widehat\varepsilon(uE'_{(1)}E_{(1)}\di v)\,\, (S_C(E_{(2)})\di E'_{(2)})\\
&=\varphi_B(uE'_{(1)}E_{(1)})\varphi_C(v)\,\, (S_C(E_{(2)})\di E'_{(2)})\\
&=\varphi_B(E_{(1)})\varphi_C(v)\,\, (S_C(E_{(2)})uE'_{(1)}\di E'_{(2)})\\
&=\varphi_C(v)\,\, (uE'_{(1)}\di E'_{(2)}).
\end{align*}
On the other hand, if we apply the counit on the second leg of Equation (\ref{eqn:3.21}), we find
\begin{align*}
\varepsilon_s(u\di v)
&=(\iota\ot\widehat\varepsilon)((1\ot (u\di v))\widehat E)\\
&=(E'_{(1)}\di S_B(E_{(1)}))\,\, \widehat\varepsilon(u\di E_{(2)}E'_{(2)}v) \\
&=(E'_{(1)}\di S_B(E_{(1)}))\,\, \varphi_B(u)\varphi_C(E_{(2)}E'_{(2)}v) \\
&=(E'_{(1)}\di E'_{(2)}v S_B(E_{(1)}))\,\, \varphi_B(u)\varphi_C(E_{(2)}).\\
&=\varphi_B(u) E'_{(1)}\di E'_{(2)}v.
\end{align*}

We see that the ranges of the source and target maps are indeed the left and right leg of $\widehat E$, obtained as a consequence of Proposition \ref{prop:3.17}. Also observe that here, the source and target maps, have values in $\widehat A$ (and not just in $M(\widehat A)$ as in the general case). The reason for this is that actually the coproduct maps $\widehat A$ already in $\widehat A\ot\widehat A$ as we have seen before.
\snl
From Proposition 2.19 of \cite{VD-W6} we know that there are isomorphisms $\gamma_s:\varepsilon_s(A)\to \varepsilon_t(\widehat A)$ and $\gamma_t:\varepsilon_t(A)\to \varepsilon_s(\widehat A)$. For this example we get the following realizations.

\prop
Let $\gamma_s(b)=\sigma_B(b)E_{(1)}\di E_{(2)}$ for $b\in B$. Then 
\begin{equation*}
\langle cb,(u\di v)\gamma_s(b')\rangle=\langle b'cb,u\di v \rangle
\end{equation*}
for all $c,b,b'$ and $u,v$. Similarly, let $\gamma_t(c)=E_{(1)}\di E_{(2)}\sigma_C^{-1}(c)$ for all $c\in C$. Then
\begin{equation*}
\langle cb, \gamma_t(c')(u\di v)\rangle = \langle cbc',u\di v \rangle
\end{equation*}
for all $c,c',b$ and $u,v$.
\eprop

\bew
i) Let $b,u,u'\in B$ and $c,v\in C$. Then
\begin{align*}
\langle cb, (u\di v)(u'E_{(1)}\di E_{(2)})\rangle
&=\langle cb, u\di E_{(2)}\rangle \varepsilon(vu'E_{(1)})\\
&=\varphi_B(bS_C(E_{(2)}))\varphi_C(S_B(u)c)\varphi_B(S_C(v)u'E_{(1)})\\
&=\varphi_B(bS_C(E_{(2)})S_C(v)u')\varphi_C(S_B(u)c)\varphi_B(E_{(1)})\\
&=\varphi_B(bS_C(v)u')\varphi_C(S_B(u)c)\\
&=\varphi_B(\sigma_B^{-1}(u')bS_C(v))\varphi_C(S_B(u)c)\\
&=\langle c\sigma_B^{-1}(u')b,u\di v\rangle.
\end{align*}
If we replace $u'$ by $\sigma_B(u')$ we get the first result.
\snl
ii) Let $b,u\in B$  and $c,c',v\in C$. Then
\begin{align*}
\langle cb,(E_{(1)}\di E_{(2)}v')(u\di v)\rangle 
&=\langle cb,E_{(1)}\di v\rangle \varepsilon(E_{(2)}v'u)\\
&=\varphi_B(bS_C(v))\varphi_C(S_B(E_{(1)})c)\varphi_C(E_{(2)}v'S_B(u))\\
&=\varphi_B(bS_C(v))\varphi_C(v'S_B(u)S_B(E_{(1)})c)\varphi_C(E_{(2)})\\
&=\varphi_B(bS_C(v))\varphi_C(v'S_B(u)c)\\
&=\varphi_B(bS_C(v))\varphi_C(S_B(u)c\sigma_C(v'))\\
&=\langle c\sigma_C(v')b,u\di v\rangle.
\end{align*}
Now we replace $v'$ by $\sigma_C^{-1}(v')$ and we get the second result.
\ebew

In fact, it is easy to verify the following immediate consequences.

\prop\label{prop:3.20c}
The maps 
\begin{equation*}
v\mapsto E_{(1)}\di E_{(2)}v
\tussenen
u\mapsto uE_{(1)}\di E_{(2)}
\end{equation*}
are isomorphisms from the algebras $C$ and $B$ to the source algebra $\varepsilon_s(\widehat A)$ and $\varepsilon_t(\widehat A)$ respectively.
\eprop
\bew
First let $v,v'\in C$. Then
\begin{align*}
(E_{(1)}\di E_{(2)}v)(E'_{(1)}\di E'_{(2)}v')
&=\varepsilon(E_{(2)}vE'_{(1)})(E_{(1)}\di  E'_{(2)}v')\\
&=\varphi_B(S_C(E_{(2)}v)E'_{(1)})(E_{(1)}\di  E'_{(2)}v')\\
&=\varphi_B(E'_{(1)})(E_{(1)}\di  E_{(2)}vE'_{(2)}v')\\
&=E_{(1)}\di  E_{(2)}vv'.
\end{align*}
This proves the first result. 
\snl
Next let $u,u'\in B$. Similarly we have
\begin{align*}
(uE_{(1)})\di  E_{(2)})(u'E'_{(1)})\di  E'_{(2)})
&=\varepsilon(E_{(2)}u'E'_{(1)})(uE_{(1)}\di E'_{(2)})\\
&=\varphi_C(E_{(2)}S_B(u'E'_{(1)}))(uE_{(1)}\di E'_{(2)})\\
&=\varphi_C(E_{(2)})(uE_{(1)}u'E'_{(1)}\di E'_{(2)})\\
&=uu'E'_{(1)}\di E'_{(2)}. 
\end{align*}
This proves the second statement.
\ebew

Remark that we can extend the previous two results to the multiplier algebras of the source and target algebras.
\snl
A last result we verify, is the following.

\prop
If $x\in \varepsilon_t(\widehat A)$ and if $y\in \varepsilon_s(\widehat A)$ we have
\begin{equation}
(S(x)\ot 1)\widehat E)=(1\ot x)\widehat E)
\tussenen
\widehat E(y\ot 1)=\widehat E(1\ot S(y)).
\end{equation}
Here $S$ is the antipode on the dual $\widehat A$ as obtained in Proposition \ref{prop:3.10}.
\eprop

\bew
Let $x=uE_{(1)}\di E_{(2)}$ where $u\in B$. Then $S(x)=S_B^{-1}(E_{(2)})\di S_C^{-1}(uE_{(1)})$. 
From Equation (\ref{eqn:3.20}) in Proposition \ref{prop:3.17} we get
\begin{align*}
(S(x)\ot 1)\widehat E)
&=(S_B^{-1}(E''_{(2)})\di S_B(E_{(1)})) \ot  (E'_{(1)}\di E_{(2)}S_C^{-1}(uE''_{(1)})E'_{(2)}) \\
&=(E''_{(1)}\di S_B(E_{(1)})) \ot  (E'_{(1)}\di E_{(2)}E''_{(2)}S_C^{-1}(u)E'_{(2)}) \\
&=(E''_{(1)}\di S_B(E_{(1)})) \ot  (uE'_{(1)}\di E_{(2)}E''_{(2)}E'_{(2)}) \\
&=(1\ot (uE'_{(1)}\di E'_{(2)}))\widehat E.
\end{align*}
For the last equality we have used (\ref{eqn:3.21}).
\snl
Similarly, take $y=E_{(1)}\di E_{(2)}v$ where $v\in C$. Then 
\begin{equation}
S(y)=S_B^{-1}(E_{(2)}v)\di S_C^{-1}(E_{(1)})=S_B^{-1}(v)E_{(1)}\di E_{(2)}.
\end{equation}
If we insert this in Equation(\ref{eqn:3.18}) of Proposition \ref{prop:3.17} we find
\begin{align*}
\widehat E(1\ot S(y))
&=(E'_{(1)}S_B^{-1}(v)E''_{(1)}E_{(1)}\di E'_{(2)})\ot (S_C(E_{(2)})\di E''_{(2)})\\
&=(E'_{(1)}E''_{(1)}E_{(1)}\di E'_{(2)}v)\ot (S_C(E_{(2)})\di E''_{(2)})\\
&=\widehat E((E'_{(1)}\di E'_{(2)}v)\ot 1).
\end{align*}
For the last equality we have used (\ref{eqn:3.19}).
\ebew

\nl
\bf Integrals on the dual $\widehat A$ \rm
\nl
Because the intention of this paper is to illustrate the methods and the results of \cite{VD-W6}, we will use the general construction of the integrals on the dual as obtained in Proposition 3.16 of [VD-W6]. We recall the result here.

\prop\label{prop:3.21}
Take any element $a\in A$. There exists a right integral $\psi_a$ on $\widehat A$ given by
\begin{equation*}
\psi_a(\omega)=\varphi(a\varepsilon_s(c))
\qquad\quad\text{if} \quad
\omega=\varphi(\,\cdot\,c)
\end{equation*}
where $\varphi$ is a left integral on $A$ and $c\in A$.
\eprop

If we apply this result for our example, we arrive at the following property.

\prop\label{prop:3.22}
For any element $c\in C$, there is right integral $\widehat \psi$ on $\widehat A$ given by $\widehat\psi(u\di v)=\varphi_C(S_B(u)cv)$ for all $u\in B$ and $v\in C$.
\eprop

\bew
Start with an element $b\in A$. Choose any element $c_0\in C$ such that $\varphi_C(c_0)=1$. 
We now apply Proposition \ref{prop:3.21} with $a=c_0b$ and use $\widehat \psi$ for the right integral denoted in the proposition by $\psi_a$.
\snl
Consider the element $u\di v$ in $\widehat A$. For any $b'\in B$ and $c'\in C$ we have
\begin{align*}
\langle c'b',u\di v \rangle
&=\varphi_B(b'S_C(v))\varphi_C(S_B(u)c') \\
&=\varphi_B(b'S_C(v))\varphi_C(c'S_BS_CS_B(u)) \\
&=\varphi(c'b'S_BS_CS_B(u)S_C(v))
\end{align*}
where we have used that the modular automorphism of $\varphi_C$ is $S_BS_C$ (see Section \ref{s:preliminaries}) and where $\varphi$ is the left integral on $A$ defined by $\varphi(c''b'')=\varphi_C(c'')\varphi_B(b''))$ as in Proposition \ref{prop:2.2} of Section \ref{s:thewmha}.
\snl
Then, by the formula in Proposition \ref{prop:3.21} above, $\widehat\psi(u\di v)$ is given by
\begin{align*}
\widehat\psi(u\di v)
&=\varphi(c_0b\varepsilon_s(S_BS_CS_B(u)S_C(v)))\\
&=\varphi(c_0bS_CS_BS_CS_B(u)S_C(v))
\end{align*}
where we have used the formula for $\varepsilon_s$ on $A$ as recalled in Section \ref{s:thewmha}. Now we apply again the formula for $\varphi$ on $A$ and we find
\begin{equation*}
\widehat\psi(u\di v)
=\varphi_C(c_0)\varphi_B(bS_CS_BS_CS_B(u)S_C(v)).
\end{equation*}
Next we use that $\varphi_C(c_0)=1$ and we define $c$ by $b=S_CS_BS_C(c)$. Then 
\begin{align*}
\widehat\psi(u\di v)
&=\varphi_B(S_CS_BS_CS_B(S_B^{-1}(c)u)S_C(v)) \\
&=\varphi_C(vS_BS_CS_B(S_B^{-1}(c)u))\\
&=\varphi_C(S_B(S_B^{-1}(c)u)v)\\
&=\varphi_C(S_B(u)cv)
\end{align*}
This completes the proof.
\ebew
In what follows, we will denote the right integral obtained above from the element $c\in C$ by $\widehat\psi_c$.
\snl
We verify that this functional is right invariant.
\snl
Let $c\in C$ and $\widehat\psi_c$ as in the above proposition. Then for $u\in B$ and $v\in C$ we have
\begin{align*}
(\widehat\psi_c\ot\iota)\Delta(u\di v)
&=\sum_{(u),(v)}\widehat\psi_c(u_{(1)}\di v_{(1)}) u_{(2)}\di v_{(2)}\\
&=\sum_{(u),(v)}\varphi_C(S_B(u_{(1)})cv_{(1)}) u_{(2)}\di v_{(2)}\\
&=\sum_{(u)}\varphi_C(S_B(u_{(1)})cS_B(E_{(1)})) u_{(2)}\di E_{(2)}v\\
&=\sum_{(u)}\varphi_C(S_B(u_{(1)})S_B(E_{(1)})) u_{(2)}\di E_{(2)}cv\\
&=\sum_{(u)}\varphi_B(E_{(1)}u_{(1)}) u_{(2)}\di E_{(2)}cv\\
&=\sum_{(u)} u_{(2)}\di S_B(u_{(1)})cv\\
&= S_C(E_{(2)})\di S_B(uE_{(1)})cv\\
&= E_{(1)}\di E_{(2)}S_B(u)cv.
\end{align*}
This gives an element in $\widehat A_s$ proving that $\widehat\psi_c$ is right invariant.
\snl
One can verify that the above argument is still valid if $c$ is replaced by a multiplier in $x\in M(C)$. So also $\widehat\psi_x$ defined by $\widehat\psi_x(u\di v)=\varphi_C(S_B(u)xv)$ will still be a right integral. 
\snl
In the following proposition we show that any right integral is of this form for some $x\in M(C)$.

\prop\label{prop:3.23}
If $\widehat\psi$ is a right integral on $\widehat A$, 
then there is an element $x\in M(C)$ such that $\widehat\psi(u\di v)=\varphi_C(S_B(u)xv)$ for all $u\in B$ and $v\in C$.
\eprop

\bew
Assume that $\widehat\psi$ is a right integral. Then $(\widehat\psi\ot\iota)\Delta(u\di v)\in \widehat A_s$. This means that for all $u\in B$ and $v\in C$ there is an element $\gamma(u,v)\in C$ so that 
\begin{equation*}
(\widehat\psi\ot\iota)\Delta(u\di v)=E_{(1)}\di  E_{(2)}\gamma(u,v).
\end{equation*}
If we insert the formula for $\Delta(u\di v)$ we find
\begin{equation}
\widehat\psi(E'_{(1)}\di S_B(E_{(1)}))S_C(E'_{(2)})u\di E_{(2)}v=E_{(1)}\di E_{(2)}\gamma(u,v).\label{eqn:3.19b}
\end{equation}
If we apply $\varphi_B$ on the first leg of this equation we get
\begin{align}
\gamma(u,v)&=\widehat\psi(E'_{(1)}\di S_B(E_{(1)}))\varphi_B(S_C(E'_{(2)})u)E_{(2)}v \nonumber\\
&=\widehat\psi(u\di S_B(E_{(1)}))E_{(2)}v.\label{eqn:3.20b}
\end{align}
It follows that $\gamma(u,vv')=\gamma(u,v)v'$.
\snl
If on the other hand, we apply $\varphi_C$ on the second leg of (\ref{eqn:3.19b}) we get
\begin{equation*}
S_B^{-1}\gamma(u,v)=\widehat\psi(E'_{(1)}\di v)\,S_C(E'_{(2)})u.
\end{equation*}
From this it follows that $\gamma(uu',v)=S_B(u')\gamma(u,v)$.
\snl
The two results together imply that there is a multiplier $x\in M(C)$ so that $\gamma(u,v)=S_B(u)xv$. If we insert this in (\ref{eqn:3.20b}) and apply $\varphi_C$ we arrive at
\begin{equation*}
\varphi_C(S_B(u)xv)=\widehat\psi(u\di S_B(E_{(1)}))\,\varphi_C(E_{(2)}v)=\widehat\psi(u\di v).
\end{equation*}
This proves the result.
\ebew

We can take $x=1$ and then we get a faithful right integral $\widehat\psi$ as we show next.

\prop\label{prop:3.24}
Define $\widehat\psi$ on $\widehat A$ by $\widehat\psi(u\di v)=\varphi_C(S_B(u)v)$. Then $\widehat\psi$ is a faithful right integral. Its modular automorphism $\sigma'$ is given by $\sigma'=S^{-2}$.
\eprop

\bew
i) We know already from  previous results that $\widehat\psi$ is right invariant.
\snl
ii) We show that it is faithful.
Take an element $\sum_i u_i\di v_i$ in $\widehat A$. Then
\begin{equation*}
(\sum_i u_i\di v_i)(u\di v)=\sum_i \varepsilon(v_iu)u_i\di v
\end{equation*}
for all $u,v$. Now assume that $\widehat\psi$ is $0$ on all such elements. By the faithfulness of $\varepsilon$ we get that
\begin{equation*}
\sum_i\varphi_C(S_B(u_i)v)v_i=0
\end{equation*}
for all $v$. From the faithfulness of $\varphi_C$ we get that $\sum_iu_i\di v_i=0$.
\snl
Similarly
\begin{equation*}
(u\di v)(\sum_i u_i\di v_i)=\sum_i \varepsilon(vu_i)u\di v_i
\end{equation*}
for all $u,v$. If we assume that $\widehat\psi$ is $0$ on all such elements, from the faithfulness of $\varepsilon$ we now get
\begin{equation*}
\sum_i\varphi_C(S_B(u)v_i)u_i=0
\end{equation*}
for all $u$. Again, because $\varphi_C$ is faithful, we find $\sum_iu_i\di v_i=0$.
\snl
iii) For all $u,u'\in B$ and $v,v'\in C$ we have
\begin{align*}
\widehat\psi((u\di v)(u'\di v'))
&=\varepsilon(vu')\widehat\psi(u\di v')\\
&=\varphi_C(vS_B(u'))\varphi_C(S_B(u)v')\\
&=\varphi_C(S_B(u')\sigma_C(v))\varphi_C(v'\sigma_C(S_B(u)))\\
&=\varphi_C(S_B(u')\sigma_C(v))\varphi_C(v'S_B(\sigma_B^{-1}(u)))\\
&=\varepsilon(v'\sigma_B^{-1}(u))\widehat\psi(u'\di \sigma_C(v))\\
&=\widehat\psi((u'\di v')(\sigma_B^{-1}(u)\di \sigma_C(v))).
\end{align*}
This proves that $\sigma'(u\di v)=\sigma_B^{-1}(u)\di \sigma_C(v)$.
\snl
From the formulas for $\sigma_B$ and $\sigma_C$ we get 
\begin{equation*}
\sigma_B^{-1}(u)\di \sigma_C(v)=S_CS_B(u)\di S_BS_C(v)=S^{-2}(u\di v)
\end{equation*} 
because $S(u\di v)=S_B^{-1}(v)\di S_C^{-1}(u)$ by Proposition \ref{prop:3.10}. 
\snl
This completes the proof.
\ebew

We need to make some important remarks here.

\opm
i) From the general theory, we know that any right integral has to be of the form $\widehat\psi(\,\cdot\,z)$ for some $z\in \widehat A_t$.  See Section 1 in \cite{VD-W6}.
\snl 
In Proposition \ref{prop:3.20c}, we have seen that $u\mapsto uE_{(1)}\di E_{(2)}$ is an isomorphism from $B$ to the target algebra of the dual. It follows that any element in $\widehat A_t$ has the form $yE_{(1)}\di E_{(2)}$ for $y\in M(B)$. Then all elements are also of the form $E_{(1)}\di xE_{(2)}$ for some $x\in M(C))$. 
\snl
Now, take the element $x\in M(C)$ so that $z=E_{(1)}\di xE_{(2)}$. Then, for every $u,v$ we find
\begin{align*}
\widehat\psi((u\di v)z)
&=\widehat\psi(u\di xE_{(2)})\varepsilon(vE_{(1)})\\
&=\widehat\psi(u\di xE_{(2)})\varphi_B(S_C(v)E_{(1)})\\
&=\widehat\psi(u\di xv)=\varphi_C(S_B(u)xv).
\end{align*}
We recover the result of Proposition \ref{prop:3.23}.
\snl
ii) We did not get the faithful right integral $\widehat\psi$ from the general construction in Proposition \ref{prop:3.22}. Still, a similar argument as in the proof of Proposition  \ref{prop:3.24} will give that the set of integrals $\widehat\psi_c$ where $c$ varies over all the elements in $C$, is a faithful set of integrals. We know from the general theory that this should be the case.
\snl
iii) Finally, we observe the following. Consider a linear functional on $\widehat A$ of the form $u\di v\mapsto \widehat\psi((u\di v)(u'\di v'))$ for a given pair of elements $u',v'$. We get
\begin{equation*}
\widehat\psi((u\di v)(u'\di v'))=\varphi_C(S_B(u)v')\varepsilon(vu').
\end{equation*}
We know that $C$ has local units and so there is an element $c$ satisfying $cv'=v'$. This will imply that 
\begin{equation*}
\widehat\psi((u\di v)(u'\di v'))=\varphi_C(S_B(u)cv')\varepsilon(vu')=\widehat\psi_c((u\di v)(u'\di v'))
\end{equation*}
for all $u,v$. Consequently, if we look at the dual of $\widehat A$ again, constructed with the faithful set of right integrals, we get the same functionals as when we would start with the single faithful one.
\eopm

We now finish this item with the formulas for the left integrals and the relations between left and right integrals.

\prop\label{prop:3.27}
For all $y\in M(B)$ there is a left integral $\widehat\varphi_y$ on $\widehat A$ given by $\widehat\varphi_y(u\di v)=\varphi_B(uyS_C(v))$ for all $u,v$. Any left integral is of this form. When $y=1$ we have a faithful left integral $\widehat\varphi$ whose modular automorphism $\sigma$ is given by $S^2$.
\eprop

\bew
We obtain these results by using the corresponding results for the right integrals and by composing with the antipode.
\snl
If e.g.\ $x\in M(C)$ we have the right integral $\widehat\psi_x$ given by $\widehat\psi_x(u\di v)=\varphi_C(S_B(u)xv)$. Then 
\begin{equation*}
\widehat\psi_x(S(u\di v))=\varphi_C(S_B(S_B^{-1}(v))xS_C^{-1}(u))=\varphi_B(uS_C(x)S_C(v))
\end{equation*}
and with $y=S_C(x)$ we get the formula for a left integral. When $x=1$, we get a faithful integral.
\snl
The modular automorphism of $\widehat\psi$ is given by $S^{-2}$. Then the modular automorphism of $\widehat\psi\circ S$ will be given by $S^2$.
\ebew

Similar remarks can be made here for left integrals as we have before for right integrals.
\snl
And finally, we obtain the relation between the left integral $\widehat\varphi$ from Proposition \ref{prop:3.24} and the right integral $\widehat\psi$ from Proposition \ref{prop:3.27}.

\prop 
If $\widehat\varphi$ is the left integral given by $\widehat\varphi(u\di v)=\varphi_B(uS_C(v))$, then \\ $\widehat\varphi(S(u\di v))=\widehat\varphi((u\di v)\delta))$ where $\delta$ is the multiplier of $\widehat A$ given by 
\begin{equation*}
\delta(u\di v)=\sigma_B^{-2}(u)\di v
\tussenen
(u\di v)\delta=u\di \sigma_C^{-2}(v)
\end{equation*}
\eprop

\bew
For all $u,v$ we have
\begin{align*}
\widehat\varphi(S(u\di v))
&=\widehat\varphi(S_B^{-1}(v)\di S_C^{-1}(u))
=\varphi_B((S_B^{-1}(v)u)\\
&=\varphi_B(u\sigma_B(S_B^{-1}(v)))
=\varphi_B(uS_C(v'))
\end{align*}
where $v'=S_C^{-1}\sigma_BS_B^{-1}(v)= \sigma_C^{-2}(v)$.
\snl
We now verify that the maps $\sigma_B$ on $B$ and $\sigma_C$ on $C$ are adjoint of each other for the bilinear form $(u,v)\mapsto \varepsilon(vu)$. Indeed, for all $u,v$ we find 
\begin{equation*}
\varepsilon(\sigma_C(v)u)=\varphi_C(S_BS_C(v)S_B(u))=\varphi_B(uS_C(v))
\end{equation*}
on the one hand, while on the other hand also
\begin{equation*}
\varepsilon(v\sigma_B(u))=\varphi_B(S_C(v)\sigma_B(u))=\varphi_B(uS_C(v)).
\end{equation*}
This completes the proof.
\ebew

We see that 
\begin{equation*}
\delta^{-1}(u\di v)\delta=\sigma_B^2(u)\di \sigma_C^{-2}(v)=S^4(u\di v)
\end{equation*}
for all $u,v$. This looks like a special case of Radford's formula for the fourth power of the antipode on $\widehat A$ (as it is proven for algebraic quantum groups in \cite{D-VD-W}). Indeed, for this example, $A$ is unimodular in the sense that it has a single faithful left integral that is also a right integral (see Proposition \ref{prop:2.2}) so that the modular element for $A$ is trivial. 
\snl
There is however the question whether or not, there is such a result as Radford's formula for these algebraic quantum groupoids. We refer to the next section for a further reflection on this problem.

\nl

\section{\hspace{-17pt}. Conclusion and further research}\label{s:conclusions}  

In this paper, we have started with a regular separability idempotent $E$ in the multiplier algebra $M(B\ot C)$ where $B$ and $C$ are non-degenerate algebras. We have considered the canonical weak multiplier Hopf algebra $(A,\Delta)$ associated with this separability idempotent. The most elaborated part of the paper is the construction of the dual $(\widehat A,\widehat\Delta)$ and the  study of its properties. This is done in great detail with the intention of illustrating the various definitions and results of the general duality theory for weak multiplier Hopf algebras as studied in \cite{VD-W6}.
\snl
For this special case, the underlying algebra $\widehat A$ is an infinite matrix algebra, only dependent on the vector space structures of  $B$ and $C$, together with a specific pairing between the two. The complexity of the example therefore lies in the coproduct. It depends highly on the separability idempotent. Moreover, various (non-trivial) features of the theory of a weak multiplier Hopf algebra with integrals are well documented using this example. For this reason, it is also instructive and helps to understand different aspects of the general theory.
\snl
There is some obvious more research to be done.
\snl
We have not looked at the involutive structure here. However, if the algebras $B$ and $C$ are $^*$-algebras and if $E$ is a self-adjoint idempotent, we know that the corresponding weak multiplier Hopf algebra is a weak multiplier Hopf $^*$-algebra with positive integrals. One has to investigate the associated involutive structure on the dual and it is expected that the integrals on the dual can also be chosen to be positive. This allows one to move to the Hilbert space level and obtain some examples of measure quantum groupoids as in \cite{T2}.
\snl
There is also the more general construction in \cite{VD-W5} where not only a separability idempotent is considered but also a regular weak multiplier Hopf algebra. This case is considered in the more recent work \cite{T2}.
\snl
It is still not clear if there exists non-regular separability idempotents. We feel that this dual example might give some inspiration to construct such examples.

\end{document}